\documentclass[12pt,a4paper,twoside,leqno]{article}

\usepackage{amssymb}

\newcommand\qed{\hfill$\square$}
\newcommand\RR{\mathbb{R}}
\newtheorem{theo}{Theorem}[section]
\newtheorem{prop}{Proposition}[section]
\newtheorem{lem}{Lemma}[section]
\newtheorem{cor}{Corollary}[section]

\newcommand{\beqn}{\begin{equation}}
\newcommand{\eeqn}{\end{equation}}
\newcommand{\bear}{\begin{eqnarray}}
\newcommand{\eear}{\end{eqnarray}}
\newcommand{\bean}{\begin{eqnarray*}}
\newcommand{\eean}{\end{eqnarray*}}

\def\qed{\hfill{{$\square$}}}
\newcommand{\rf}[1]{(\ref{#1})}
\newcommand{\proof}{\noindent {\sc Proof.} }

{\medskip \noindent {\bf QUESTION:~}%
\sf }%
{\medskip}

\newcounter{remar}
\newenvironment{remark}%
{\stepcounter{remar}\medskip \noindent {\it Remark \thesection.\arabic{remar}}%
\hspace{1.5mm}}%
{\par}

\title{\Large\bf Asymptotic profiles of solutions \\
to viscous Hamilton-Jacobi equations}

\author{
Sa\"\i d Benachour\\
{\small Institut Elie Cartan-Nancy, Universit\'e Henri Poincar\'e }\\
{\small BP 239, F-54506 Vand\oe uvre l\`es Nancy cedex, France} \\
{\small E-mail: {\it benachou@iecn.u-nancy.fr}} \\
\\
Grzegorz Karch\\
{\small Instytut Matematyczny, Uniwersytet Wroc\l awski }\\
{\small pl. Grunwaldzki 2/4, 50-384 Wroc\l aw, Poland }\\
{\small Institute of Mathematics, Polish Academy} \\ 
{\small of Sciences, Warsaw (2002-2003)} \\
{\small E-mail: {\it karch@math.uni.wroc.pl}}\\
\\ 
Philippe Lauren\c cot\\
{\small Math\'ematiques pour l'Industrie et la Physique} \\
{\small CNRS UMR 5640, Universit\'e Paul Sabatier-Toulouse 3}\\
{\small 118 route de Narbonne, F-31062 Toulouse cedex 4, France} \\
{\small E-mail: {\it laurenco@mip.ups-tlse.fr}}
} 

\begin{document}
\maketitle

\begin{abstract}

The large time behavior of solutions to the Cauchy problem for the
viscous Hamilton-Jacobi equation $u_t-\Delta u+|\nabla u|^q=0$ is
classified. If $q>q_c:= (N+2)/(N+1)$, it is shown that non-negative
solutions corresponding to integrable initial data converge in
$W^{1,p}(\RR^N)$ as $t\to\infty$ toward a multiple of the fundamental
solution for the heat equation for every $p\in
[1,\infty]$ (diffusion-dominated case). On the other hand, if
$1<q<q_c$, the large time asymptotics is given by the very singular
self-similar solutions of the viscous Hamilton-Jacobi equation. 

For non-positive and integrable solutions, the large time behavior of
solutions is more complex. The case $q\geq 2$ corresponds to the
diffusion-dominated case. The diffusion profiles in the large time
asymptotics appear also for $q_c<q<2$ provided suitable smallness
assumptions are imposed on the initial data. Here, however, the most
important result asserts that under some conditions on initial
conditions and for $1<q<2$, the large time behavior of solutions is
given by the self-similar viscosity solutions to the non-viscous
Hamilton-Jacobi equation $z_t+|\nabla z|^q=0$ supplemented with the
initial datum $z(x,0)=0$ if $x\neq 0$ and $z(0,0)<0$.   

\end{abstract}

\vspace{0.5cm}

\centerline{R\'esum\'e}

\small

Nous classifions le comportement asymptotique des solutions du 
probl\`eme de Cauchy pour l'\'equation de Hamilton-Jacobi avec 
diffusion $u_t-\Delta u+|\nabla u|^q=0$. Si $q>q_c:= (N+2)/(N+1)$, 
nous montrons que, lorsque $t\to\infty$, les solutions int\'egrables 
et positives convergent dans $W^{1,p}(\RR^N)$ vers un multiple de la 
solution fondamentale de l'\'equation de la chaleur pour tout $p\in
[1,\infty]$ (diffusion dominante). Ensuite, si $1<q<q_c$, le comportement 
asymptotique est d\'ecrit par la solution tr\`es singuli\`ere auto-similaire 
de l'\'equation de Hamilton-Jacobi avec diffusion. 

En ce qui concerne les solutions int\'egrables et n\'egatives, la 
situation est plus complexe. Le terme de diffusion est de nouveau 
dominant si $q\geq 2$, ainsi que lorsque $q_c<q<2$ pourvu que la donn\'ee 
initiale soit suffisamment petite. Ensuite, pour $1<q<2$, nous identifions 
une classe de donn\'ees initiales pour laquelle le comportement 
asymptotique des 
solutions est donn\'e par une solution de viscosit\'e auto-similaire 
de l'\'equation de Hamilton-Jacobi $z_t+|\nabla z|^q=0$ avec la condition 
initiale (non continue) $z(x,0)=0$ si $x\neq 0$ et $z(0,0)<0$. 

\normalsize

\vspace{0.2cm}  


\vspace{0.2cm}  

\noindent Mots-cl\'es~: Equation de Hamilton-Jacobi diffusive, 
comportement asymptotique auto-similaire, estimations 
unilat\'erales du Laplacien.

\vspace{0.5cm}  

\section{Introduction}

We investigate the large time behavior of integrable solutions to the
Cauchy problem for the viscous Hamilton-Jacobi equation  
\bear
&&u_t-\Delta u +|\nabla u|^q=0\,, \quad x\in \RR^N\,, \quad
t>0\,,\label{eq}\\ 
&&u(x,0)=u_0(x)\,, \quad x\in\RR^N\,,\label{ini}
\eear
where $q>1$. The dynamics of the solutions to (\ref{eq})-(\ref{ini})
is governed by two competing effects, namely those resulting from the
diffusive term $-\Delta u$ and those corresponding to the
``hyperbolic'' nonlinearity $\vert\nabla u\vert^q$. Our aim here is to
figure out whether one of these two effects rules the large time
behavior, according to the values of $q$ and the initial data
$u_0$. Since the nonlinear term $\vert\nabla u\vert^q$ is
non-negative, it acts as an absorption term for non-negative solutions
and as a source term for non-positive solutions. We thus consider
separately non-negative and non-positive solutions. Let us outline our
main results now.  

For non-negative initial data, it is already known that diffusion
dominates the large time behavior for $q>q_c:=(N+2)/(N+1)$ and that
the nonlinear term only becomes effective for $q<q_c$
\cite{AB98,BL99,BK99,BGK03}. We obtain more precise information in
Theorems \ref{th:01} and  \ref{th:02} below. In particular, if $q\in
(1,q_c)$ and the initial datum decays sufficiently rapidly at
infinity, there is a balance between the diffusive and hyperbolic
effects: the solution $u(t)$ behaves for large $t$ like the very
singular solution to (\ref{eq}), the existence and uniqueness of which
have been established in \cite{BL01,BKLxx,QW01}.  

For non-positive initial data, there are two critical exponents
$q=q_c$ and $q=2$, as already noticed in \cite{LS03}, and the picture
is more complicated. More precisely, the diffusion governs the large
time dynamics for any initial data if $q\ge 2$ and for sufficiently
small initial data if $q\in (q_c,2)$, and we extend the result from
\cite[Proposition 2.2]{LS03} in that case (cf. Theorem \ref{th:03},
below). On the other hand, when $q\in (1,2)$, we prove that, for
sufficiently large initial data, the large time behavior is governed
by the nonlinear reaction term. This fact  is also true for any
initial datum $u_0\not\equiv 0$ if $N\le 3$ and $q$ is sufficiently
close to $1$. We actually conjecture that the nonlinear reaction term
always dominates in the large time for any non-zero initial datum as
soon as $q\in (1,q_c)$.  

Let us finally mention that, when $q\in (q_c,2)$, there is at least
one (self-similar) solution for which there is a balance between the
diffusive and hyperbolic effects for large times \cite{BSW02}.  

Before stating more precisely our results, let us recall that for
every initial datum $u_0\in W^{1,\infty}(\RR^N)$ the Cauchy problem
(\ref{eq})-(\ref{ini}) has a unique global-in-time solution which is
classical for positive times, that is  
$$
u\in\mathcal{C}(\RR^N\times [0,\infty)) \cap
\mathcal{C}^{2,1}(\RR^N\times (0,\infty))\,. 
$$ 
In addition, this solution satisfies the estimates
\beqn
\label{contrex}
\|u(t)\|_\infty\le \|u_0\|_\infty \quad \mbox{and}\quad 
\|\nabla u(t)\|_\infty\le \|\nabla u_0\|_\infty \quad \mbox{for all
$t>0$.}
\eeqn
Moreover, by the maximum principle, $u_0\ge 0$ implies that $u\ge 0$
and $u_0\le 0$ ensures that $u\le 0$.   
We refer the reader to \cite{AB98,BL99,GGK03} for the proofs of all
these preliminary results. In addition, a detailed analysis of the
well-posedness of (\ref{eq})-(\ref{ini}) in the Lebesgue spaces
$L^p(\RR^N)$ may be found in the recent paper \cite{BSW02}. 

\bigskip

{\bf Notations.}
The notation to be used is mostly standard. For $1\le p\le \infty$,
the $L^p$-norm of a Lebesgue measurable real-valued function $v$
defined on $\RR^N$ is denoted by $\|v\|_p$. We will always denote by
$\|\cdot\|_{\cal X}$ the norm of any other Banach space $\cal X$ used
in this paper. Also, $W^{1,\infty}(\RR^N)$ denotes the Sobolev space
consisting of functions in $L^\infty (\RR^N)$ whose first order
generalized derivatives belong  to $L^\infty(\RR^N)$. The space of
compactly supported and $\mathcal{C}^\infty$-smooth functions in
$\RR^N$ is denoted by $\mathcal{C}_c^\infty(\RR^N)$, and
$\mathcal{C}_0(\RR^N)$ is the set of continuous functions $u$ such
that 
$$
\lim_{R\to\infty} \sup_{\vert x\vert\ge R}\, \{\vert u(x)\vert\} = 0\,.
$$
For a real number $r$, we denote by $r^+:=\max{\{r,0\}}$ its positive
part and by $r^-:=\max{\{-r,0\}}$ its negative part. The letter $C$
will denote generic positive constants, which do not depend on $t$ and
may vary from line to line during computations. Throughout the paper,
we use the critical exponent  
$$
q_c := \frac{N+2}{N+1}\,.
$$

%%%%%%%%%%%%%%%%%%%%%%%%%%%%%%%%%%%%%%%%%%%%%%%%%%%%%%%%%%%%%%%%%%%%%%%%%%
%%%%%%
%%%%%%%%%%%%%%%%%%%%%%%%%%%%%%%%%%%%%%%%%%%%%%%%%%%%%%%%%%%%%%%%%%%%%%%%%%
%%%%%%

\section{Results and comments}
\setcounter{equation}{0}
\setcounter{remar}{0}

As already outlined, the large time behavior of solutions to
\rf{eq}-\rf{ini} is determined not only by the exponent $q$ of the
nonlinear term $\vert\nabla u\vert^q$ but also by the sign, size, and
shape of the initial conditions. In the present paper, we attempt to
describe this variety of different asymptotics of solutions, imposing
particular assumptions on initial data. In order to present our
results in the most transparent form, we divide this section into
subsections. 

\subsection{Non-negative initial conditions}

In Theorems \ref{th:01} and \ref{th:02} below, we always assume that 
\begin{equation}\label{u0+}
u_0 \;\mbox{ is a non-negative function in }\; L^1(\RR^N)\cap
W^{1,\infty}(\RR^N)\,, \quad u_0\not\equiv 0\,, 
\end{equation}
and we denote by $u=u(x,t)$ the corresponding non-negative solution of
the Cauchy problem \rf{eq}-\rf{ini}. In that case, we recall that
$t\longmapsto \Vert u(t)\Vert_1$ is a non-increasing function and that
$\vert\nabla u\vert$ belongs to $L^q(\RR^N\times (0,\infty))$. In
addition,  
\beqn
I_\infty := \lim_{t\to\infty}\int_{\RR^N}u(x,t)\;dx=
\int_{\RR^N}u_0(x)\;dx-\int_0^\infty \int_{\RR^N} |\nabla
u(x,s)|^q\;dx\,ds
\label{I:infty}
\eeqn
satisfies $I_\infty>0$ if $q>q_c$ and $I_\infty=0$ if $q\le q_c$
(cf. \cite{AB98,BL99,BK99}, for details). Since we would have
$I_\infty=\Vert u_0\Vert_1>0$ for the linear heat equation, we thus
say that diffusion dominates the large time behavior when
$I_\infty>0$, that is, when $q>q_c$.  

We first consider the diffusion-dominated case.

\begin{theo} \label{th:01}
Suppose \rf{u0+} and that $q>q_c$. For every $p\in [1,\infty]$,
\beqn
\lim_{t\to\infty} t^{(N/2)(1-1/p)} \|u(t)-I_\infty
G(t)\|_p=0\label{lim:G1}
\eeqn
and
\beqn
\lim_{t\to\infty} t^{(N/2)(1-1/p)+1/2} \|\nabla u(t)-I_\infty \nabla
G(t)\|_p=0.\label{lim:G2}
\eeqn
Here, $G(x,t)=(4\pi t)^{-N/2} \exp(-|x|^2/(4t))$ is the fundamental
solution of the heat equation.
\end{theo}

When $p=1$, the relation \rf{lim:G1} is proved in \cite{BGK03} and
Theorem \ref{th:01} extends the convergence of $u(t)$ towards a
multiple of $G(t)$ to $W^{1,p}(\RR^N)$, $p\in [1,\infty]$.  

\begin{remark}
Theorem \ref{th:01} holds true when $I_\infty=0$ (i.e. for $q\le q_c$)
as well, but in that case, the relation \rf{lim:G1} says only that
$\|u(t)\|_p$ tends to 0 as $t\to\infty$ faster than
$t^{-(N/2)(1-1/p)}$. 
\end{remark}

\bigskip
Our next theorem is devoted to the balance case $1<q<q_c$  when a
particular self-similar solution of \rf{eq} appears in the large time
asymptotics. 

\begin{theo}\label{th:02}
Suppose \rf{u0+}.
Assume that $q\in (1,q_c)$ and, moreover, that 
\begin{equation}
\mathop{\mbox{\rm ess lim}}_{|x|\to\infty} |x|^a u_0(x)=0 \quad
\mbox{with}\quad a={2-q\over q-1}. \label{a-u0} 
\end{equation}
For every $p\in [1,\infty]$,
\beqn
\lim_{t\to\infty} t^{(N/2)(1-1/p) +(a-N)/2} \|u(t)- W(t)\|_p=0
\label{lim:W1}
\eeqn
and
\beqn
\lim_{t\to\infty} t^{(N/2)(1-1/p)+(a-N)/2+1/2} \|\nabla u(t)-
\nabla W(t)\|_p=0,
\label{lim:W2}
\eeqn
where $W(x,t)=t^{-a/2} W(x t^{-1/2},1)$ is the very singular self-similar
solution to \rf{eq}.
\end{theo}

For the existence and uniqueness of the very singular solution to
\rf{eq}, we refer the reader to \cite{BL01,BKLxx,QW01}. 
Notice also that the initial datum $u_0$ is integrable by assumption
\rf{a-u0} since $a>N$ for $1<q<q_c$. 

\begin{remark}
In the critical case $q=q_c$, it is also expected that $u(t)$
converges towards a multiple of $G(t)$ with a correction in the form
of an extra logarithmic factor resulting from the absorption
term. This conjecture is supported by what is already known for
non-negative solutions to the Cauchy problem $w_t - \Delta w +
w^{(N+2)/N}=0$ (see, e.g., \cite{Va93} and the references therein).  
\end{remark}

%%%%%%%%%%%%%%%%%%%%%%%%%%%%%%%%%%%%%%%%%%%%%%%%%%%%%%%%%%%%%%%%%%%%%%%%%%
%%%%%%
%%%%%%%%%%%%%%%%%%%%%%%%%%%%%%%%%%%%%%%%%%%%%%%%%%%%%%%%%%%%%%%%%%%%%%%%%%
%%%%%%

\subsection{Non-positive initial conditions}

We now turn to non-positive solutions and assume that 
\begin{equation}\label{u0-}
u_0 \;\mbox{ is a non-positive function in }\; L^1(\RR^N)\cap
W^{1,\infty}(\RR^N)\,, \quad u_0\not\equiv 0\,. 
\end{equation}
We denote by $u=u(x,t)$ the corresponding non-positive solution of the
Cauchy problem \rf{eq}-\rf{ini}. In that case, we recall that
$t\longmapsto \Vert u(t)\Vert_1$ is a non-decreasing function and put 
\beqn
I_\infty := \inf_{t\ge 0}\int_{\RR^N}u(x,t)\;dx = - \sup_{t\ge 0}\Vert
u(t)\Vert_1 \in [-\infty,-\Vert u_0\Vert_1]\,. 
\label{evian}
\eeqn

Substituting $u=-v$ in  \rf{eq}-\rf{ini} we obtain that $v=v(x,t)$ is
a non-negative solution to  
\begin{equation}
v_t-\Delta v-|\nabla v|^q=0, \quad v(x,0)= - u_0(x),\label{eq-}
\end{equation} 
which has been studied in \cite{BSW02, GGKpp, GGK03, LS03}.

\bigskip

We start again with the diffusion-dominated case.

\begin{theo}\label{th:03}
Suppose \rf{u0-}.

a) Assume that $q\ge 2$. Then $I_\infty>-\infty$ and $\vert\nabla
u\vert$ belongs to $L^q(\RR^N\times (0,\infty))$. In addition,
$I_\infty$ is given by \rf{I:infty} and the relations \rf{lim:G1} and
\rf{lim:G2} hold true for every $p\in [1,\infty]$.  

b) Assume that $q\in (q_c,2)$. There exists
$\varepsilon=\varepsilon(N,q)$ such that, if   
\beqn
\|u_0\|_1\|\nabla u_0\|_\infty ^{(N+1)q-(N+2)}<\varepsilon\,,
\label{small:u0}
\eeqn
then the conclusions of part a) are still valid. 
\end{theo}

The fact that $I_\infty>-\infty$ under the assumptions of Theorem
\ref{th:03} is established in \cite{LS03}, together with the relation
\rf{lim:G1} for $p=1$. We extend here this convergence to
$W^{1,p}(\RR^N)$, $p\in [1,\infty]$. 

The smallness assumption imposed in  \rf{small:u0} is necessary to
obtain the heat kernel as the first term of the asymptotic expansion
of solutions. This is an immediate consequence of the following
theorem and the subsequent discussion. 

\begin{theo}\label{th:04}
Suppose \rf{u0-} and that $q\in (q_c,2)$.

a) There exists a non-positive self-similar solution 
$$
V=V(x,t)=t^{-(2-q)/(2(q-1))}V(xt^{-1/2},1)
$$
 to \rf{eq} such that
$$
\lim_{t\to\infty} t^{(N/2)(1-1/p)}\|V(t)\|_p =\infty \quad\mbox{ and
}\quad \lim_{t\to\infty} t^{(N/2)(1-1/p)+1/2}\|\nabla V(t)\|_p =
\infty 
$$
for all $p\in [1,\infty]$.

b) There is a constant $K=K(q)\ge 0$ such that, if $u_0\in
W^{2,\infty}(\RR^N)$ satisfies
\beqn
\|u_0\|_\infty \ \left\Vert (\Delta u_0)^+
\right\Vert_{\infty}^{1-2/q} >K\label{large:u0} 
\eeqn
then
\beqn
\lim_{t\to\infty} \|u(t)\|_\infty >0.\label{lim:u:infty}
\eeqn
\end{theo}

The first part of Theorem \ref{th:04} is proved in \cite{BSW02} while
the second assertion is new. Let us point out here that, for the
Hamilton-Jacobi equation $w_t+\vert\nabla w\vert^q=0$, the
$L^\infty$-norm of solutions remains constant throughout time
evolution, while it decays to zero for the linear heat equation. We
thus realize that, under the assumptions of Theorem \ref{th:04} b),
the diffusive term is not strong enough to drive the solution to zero
in $L^\infty$ as $t\to \infty$ and the large time dynamics is
therefore ruled by the Hamilton-Jacobi term $\vert\nabla u\vert^q$. 

Unfortunately, the conditions \rf{small:u0} and \rf{large:u0} do not
involve the same quantities. Still, we can prove that if $u_0$ fulfils  
$$
\|u_0\|_\infty \ \|D^2 u_0\|_\infty^{1-2/q} >K
$$
(which clearly implies \rf{large:u0} since $q<2$), the quantity
$\|u_0\|_1\|\nabla u_0\|_\infty ^{(N+1)q-(N+2)}$ cannot be
small. Indeed, there is a constant $C$ depending only on $q$ and $N$
such that 
\begin{equation}
\left( \Vert u_0\Vert_{\infty} \ \|D^2
  u_0\|_\infty^{1-2/q}\right)^{q(N+1)/2} \le C\Vert u_0\Vert_{1} \Vert
\nabla u_0\Vert_{\infty}^{q(N+1)-(N+2)}. 
\label{2:est}
\end{equation}
For the proof of \rf{2:est}, put $B=\Vert u_0\Vert_{\infty} \|D^2
u_0\|_\infty^{1-2/q}$ and note that the Gagliardo-Nirenberg
inequalities 
\bean
\Vert u_0\Vert_{\infty} & \le & C\ \Vert \nabla
u_0\Vert_{\infty}^{N/(N+1)}\ \Vert u_0\Vert_{1}^{1/(N+1)} \,,\\
\Vert \nabla u_0\Vert_{\infty} & \le & C\ \Vert D^2
u_0\Vert_{\infty}^{(N+1)/(N+2)}\ \Vert u_0\Vert_{1}^{1/(N+2)}\,,
\eean
imply that
\bean
\Vert \nabla u_0\Vert_{\infty}^{(2-q)(N+2)} & \le & C\ \Vert D^2
u_0\Vert_{\infty}^{(2-q)(N+1)}\ \Vert u_0\Vert_{1}^{2-q} \\ 
& = & C \ B^{-q(N+1)}\ \Vert u_0\Vert_{\infty}^{q(N+1)}\ \Vert
u_0\Vert_{1}^{2-q} \\
& \le & C\   B^{-q(N+1)} \ \Vert \nabla u_0\Vert_{\infty}^{qN} \ \ \Vert
u_0\Vert_{1}^{2}\,,
\eean
whence the above claim. 

We next show that the second assertion of Theorem \ref{th:04} is also
true when $q\in (1,q_c)$.  

\begin{theo}\label{th:05}
Suppose \rf{u0-} and 
that $q\in (1,q_c]$. There is a constant $K=K(q)\ge 0$ such that, if
$u_0\in W^{2,\infty}(\RR^N)$ fulfils \rf{large:u0}, then
\rf{lim:u:infty} holds true.  

Furthermore, if $N\le 3$ and $1 < q < 4/\left( 1+\sqrt{1+2N} \right)$,
then $K(q)=0$. 
\end{theo}

We actually conjecture that $K(q)=0$ for any $q\in (1,q_c)$, but we
have yet been unable to prove it.  

\bigskip

The last result confirms the domination of the Hamilton-Jacobi term
for large times when \rf{lim:u:infty} holds true and provides precise
information on the large time behavior.

\begin{theo}\label{th:06}
Let $q\in (1,2)$. Assume that $u_0\in\mathcal{C}_0(\RR^N)$ fulfils
(\ref{u0-}) and is such that  
\beqn
\label{z1}
M_\infty := \lim_{t\to\infty} \Vert u(t)\Vert_\infty>0\,.
\eeqn
Then
\beqn
\label{z2}
\lim_{t\to\infty} \Vert u(t) - Z_{M_\infty}(t)\Vert_\infty = 0\,,
\eeqn
where $Z_{M_\infty}$ is given by 
\beqn
\label{z3}
Z_{M_\infty}(x,t) := - \left( M_\infty - (q-1)\ q^{-q/(q-1)}\ \left(
    \frac{\vert x\vert}{t^{1/q}} \right)^{q/(q-1)} \right)^+  
\eeqn
for $(x,t)\in\RR^N\times (0,\infty)$. In fact, $Z_{M_\infty}$ is the
    unique viscosity solution in $\mathcal{BUC}(\RR^N\times
    (0,\infty))$ to  
\beqn
\label{z4}
z_t+\vert\nabla z\vert^q=0 \quad \mbox{ in } \quad \RR^N\times
(0,\infty) 
\eeqn
with the bounded and lower semicontinuous initial datum $z(x,0)=0$ if
$x\ne 0$ and $z(0,0)=-M_\infty$. 
\end{theo}

\noindent The last assertion of Theorem~\ref{th:06} follows from
\cite{St02}. Moreover, $Z_{M_\infty}$ is actually given by the
Hopf-Lax formula 
$$
Z_{M_\infty}(x,t) = \inf_{y\in\RR^N} \left\{ -M_\infty\
  \mathbf{1}_{\{0\}}(y) + (q-1)\ q^{-q/(q-1)} \vert
  x-y\vert^{q/(q-1)}\ t^{-1/(q-1)} \right\} 
$$
for $(x,t)\in\RR^N\times (0,\infty)$, where $\mathbf{1}_{\{0\}}$
denotes the characteristic function of the set $\{0\}$. Observe that
$Z_{M_\infty}$ is a self-similar solution to (\ref{z4}) since
$Z_{M_\infty}(x,t)=Z_{M_\infty}(x t^{-1/q},1)$.

\medskip

If $N=1$, the convergence stated in Theorem~\ref{th:06} extends to the
gradient of $u$. 

\begin{prop}\label{pr:07}
Assume that $N=1$ and consider a non-positive function $u_0$ in
$W^{1,1}(\RR)\cap W^{1,\infty}(\RR)$. Under the assumptions and
notations of Theorem~\ref{th:06}, we have also 
$$
\lim_{t\to\infty} t^{(1-1/p)/q}\ \Vert u_x(t) -  Z_{M_\infty,x}(t)
\Vert_p = 0
$$
for $p\in[1,\infty)$.
\end{prop}

\medskip

In fact, if $N=1$ and $u_0\in W^{1,1}(\RR)\cap W^{1,\infty}(\RR)$, the
function $U:=u_x$ is a solution to the convection-diffusion equation 
\beqn
\label{cdeq}
U_t - U_{xx} + \left( \vert U\vert^q \right)_x = 0\,, \quad x\in\RR\,,
\quad t>0\,, 
\eeqn
with initial datum $U(0) = u_{0x}$ and satisfies
\beqn
\label{zeromean}
\int_{\RR} U(x,t)\ dx = \int_{\RR} u_{0x}(x)\ dx = 0\,, \quad t\ge 0\,.
\eeqn
The large time behavior of non-negative or non-positive integrable
solutions to (\ref{cdeq}) is now well-identified \cite{EZ91,EVZ1} but
this is far from being the case for solutions satisfying
(\ref{zeromean}). In this situation, some sufficient conditions on
$U(0)$ are given in \cite{KS02} for the solution to (\ref{cdeq}) to
exhibit a diffusion-dominated large time behavior. Also, convergence
to $N$-waves is studied in \cite{Kim} but, for solutions satisfying
(\ref{zeromean}), no condition is given in that paper which guarantees
that $U(t)$ really behaves as an $N$-wave for large times. As a
consequence of our analysis, we specify such a condition and also
provide several new information on the large time behavior of
solutions to (\ref{cdeq}) satisfying (\ref{zeromean}). 
Results on the large time behavior of solutions to equation \rf{cdeq}
satisfying the condition \rf{zeromean} are reviewed in the companion paper
\cite{BKLyy}. 

\medskip

We finally outline the contents of the paper: the next section is
devoted to some preliminary estimates. Theorems~\ref{th:01} and~\ref{th:03}
(diffusion-dominated case) are proved in Section~4 and Theorem~\ref{th:02} in
Section~5. The remaining sections are devoted to the
``hyperbolic''-dominated case: Theorems~\ref{th:04} and~\ref{th:05}
are proved in Section~5 and Theorem~\ref{th:06} and
Proposition~\ref{pr:07} in Section~6. 

%%%%%%%%%%%%%%%%%%%%%%%%%%%%%%%%%%%%%%%%%%%%%%%%%%%%%%%%%%%%%%%%%%%%%%%%%%
%%
%%%%%%%%%%%%%%%%%%%%%%%%%%%%%%%%%%%%%%%%%%%%%%%%%%%%%%%%%%%%%%%%%%%%%%%%%%
%%

\section{Preliminary estimates}
\setcounter{equation}{0}
\setcounter{remar}{0}

Let us first state a gradient estimate for solutions to \rf{eq} which
is a consequence of \cite[Theorem~1]{BL99} (see also \cite[Theorem~2]{GGK03}).
Note that, in this section, we do not impose a sign condition on
the solution $u$ to \rf{eq}. 

\begin{prop}\label{prestgrad}
Assume that $u=u(x,t)$ is the solution to \rf{eq}-\rf{ini}
corresponding to the initial datum $u_0 \in W^{1,\infty}(\RR^N)$. 
For every $q>1$, there is a constant $C_1>0$ depending only on $q$ such that
\beqn
\label{grad}
\Vert\nabla u(t)\Vert_{\infty} \le C_1\ \Vert u_0\Vert_{\infty}^{1/q}\
t^{-1/q}\,, \quad\mbox{for all} \quad  t>0\,.
\eeqn
\end{prop}

\proof
 Setting $v=u+\Vert u_0\Vert_{\infty}$, it
readily follows from (\ref{eq}) and the maximum principle that $v$ is a
non-negative solution to (\ref{eq}). By \cite[Theorem 1]{BL99}, there is a
constant $C$ depending only on $q$ such that 
$$
\left\Vert \nabla v^{(q-1)/q}(t)\right\Vert_{\infty} \le C\ t^{-1/q}\,,
\qquad t>0\,.
$$
Since $\nabla v = (q/(q-1))\ v^{1/q}\ \nabla v^{(q-1)/q}$ and $\vert
u(x,t)\vert\le \Vert u_0\Vert_{\infty}$, we further deduce that 
$$
\Vert \nabla u(t) \Vert_{\infty} = \Vert \nabla v(t) \Vert_{\infty}
\le C\ \Vert v(t)\Vert_{\infty}^{1/q}\ \left\Vert \nabla
v^{(q-1)/q}(t)\right\Vert_{\infty} \le C\ \Vert
u_0\Vert_{\infty}^{1/q}\ t^{-1/q}\,,
$$
whence (\ref{grad}). \hfill $\square$

\medskip

Next, we derive estimates for the second derivatives of solutions to
\rf{eq}-\rf{ini} when $q\in (1,2]$. 

\begin{prop}\label{prestlap}
Under the assumptions of Proposition \ref{prestgrad},
if $q\in (1,2]$, the Hessian matrix $D^2 u = \left( u_{x_i x_j}
\right)_{1\le i,j\le N}$ of $u$ satisfies
\bear
\label{hess1}
D^2 u(x,t) & \le & {\Vert\nabla u_0\Vert_{\infty}^{2-q}\over q\ (q-1)\
t}\ Id\,, \\
\label{hess2}
D^2 u(x,t) & \le & {C_2\ \Vert u_0\Vert_{\infty}^{(2-q)/q}\over
t^{2/q}}\ Id\,,
\eear
for $(x,t)\in\RR^N\times (0,\infty)$, where $C_2$ is a positive constant
depending only on $q$. 

Furthermore, if $u_0\in W^{2,\infty}(\RR^N)$, 
\beqn
\label{hess3}
D^2 u(x,t) \le \Vert D^2 u_0\Vert_{\infty}\ Id\,.
\eeqn
\end{prop}

In Proposition \ref{prestlap}, $Id$ denotes the identity matrix of
$\mathcal{M}_N(\RR)$. Given two matrices $A$ and $B$ in
$\mathcal{M}_N(\RR)$, we write $A\le B$ if $A\xi\cdot\xi \le
B\xi\cdot\xi$ for every vector $\xi\in \RR^N$. 

For $q=2$, the estimates \rf{hess1} and \rf{hess2} follow from the
analysis of Hamilton \cite{Ha93} (since, if $f$ is a non-negative
solution to the linear heat equation $f_t = \Delta f$, the function
$-\ln{f}$ solves (\ref{eq}) with $q=2$). In Proposition \ref{prestlap}
above, we extend that result to any $q\in (1,2]$.  
 
\begin{remark}
The estimates (\ref{hess1}) and (\ref{hess2}) may also be seen as an
extension to a multidimensional setting of a weak form of the Oleinik type
gradient estimate for scalar conservation laws. Indeed, if $N=1$ and
$U=u_x$, then $U$ is a solution to $U_t - U_{xx} + \left( \vert U \vert^q
\right)_x=0$ in $\RR\times (0,\infty)$. The estimates (\ref{hess1}) and
(\ref{hess2}) then read 
$$
U_x \le C\ \Vert U(0)\Vert_{\infty}^{2-q}\ t^{-1} \;\;\mbox{ and }\;\; U_x
\le C\ \Vert u_0\Vert_{\infty}^{(2-q)/q}\ t^{-2/q}
$$
for $t>0$, respectively, and we thus recover the results of
\cite{FL99,Kim} in that case.
\end{remark}

\medskip

\noindent
{\sc Proof of Proposition \ref{prestlap}.}  For $1\le i,j \le N$, we
put $w_{ij} = u_{x_i x_j}$. It follows from equation (\ref{eq}) that 
\begin{eqnarray}
w_{ij,t} - \Delta w_{ij} & = & -q\ \left( \vert\nabla u\vert^{q-2}\ \left(
\sum_{k=1}^N u_{x_k}\ w_{jk} \right) \right)_{x_i}\nonumber\\
& = & -q\ \vert\nabla u\vert^{q-2}\ \sum_{k=1}^N w_{ik}\ w_{jk} - q\
\vert\nabla u\vert^{q-2}\ \sum_{k=1}^N u_{x_k}\ w_{jk,x_i} \label{prop:1}\\
& & -\ q\ (q-2)\ \vert\nabla u\vert^{q-4}\ \left( \sum_{k=1}^N u_{x_k}\
w_{ik} \right)\ \left( \sum_{k=1}^N u_{x_k}\ w_{jk} \right)\,.\nonumber
\end{eqnarray}
Consider now $\xi\in\RR^N\setminus\{ 0\}$ and set 
$$ h = \sum_{i=1}^N \sum_{j=1}^N w_{ij}\ \xi_i\ \xi_j\,.$$
Multiplying \rf{prop:1} by $\xi_i\ \xi_j$ and summing up the
resulting identities yield
\begin{eqnarray}
h_t - \Delta h & = & - q\ \vert\nabla u\vert^{q-2}\ \sum_{k=1}^N \left(
\sum_{i=1}^N w_{ik}\ \xi_i \right)^2 - q\ \vert\nabla u\vert^{q-2}\ \nabla
u \cdot \nabla h \nonumber\\
& & -\ q\ (q-2)\ \vert\nabla u\vert^{q-4}\ \left( \sum_{i=1}^N
\sum_{j=1}^N u_{x_j}\ w_{ij}\ \xi_i \right)^2\,.\label{prop:2}
\end{eqnarray}
Thanks to the following inequalities
\bean
\vert\nabla u\vert^{q-4}\ \left( \sum_{i=1}^N \sum_{j=1}^N u_{x_j}\
w_{ij}\ \xi_i \right)^2 & \le & \vert\nabla u\vert^{q-4}\ \sum_{j=1}^N
\vert u_{x_j}\vert^2\ \sum_{j=1}^N \left( \sum_{i=1}^N  w_{ij}\ \xi_i
\right)^2 \\
& \le & \vert\nabla u\vert^{q-2}\ \sum_{k=1}^N \left( \sum_{i=1}^N
w_{ik}\ \xi_i \right)^2\,,
\eean
and
$$
h^2 \le \vert\xi\vert^2\ \sum_{k=1}^N \left( \sum_{i=1}^N w_{ik}\ \xi_i
\right)^2\,,
$$
and since $q\le 2$, the right-hand side of identity \rf{prop:2} can be
bounded from above. We thus obtain
\bean
h_t - \Delta h & \le & - q\ (q-1)\ \vert\nabla u\vert^{q-2}\ \sum_{k=1}^N
\left( \sum_{i=1}^N w_{ik}\ \xi_i \right)^2 - q\ \vert\nabla u\vert^{q-2}\
\nabla u \cdot \nabla h \\
& \le & - q\ \vert\nabla u\vert^{q-2}\ \nabla u \cdot \nabla h - {q\ (q-1)\
\vert\nabla u\vert^{q-2}\over \vert\xi\vert^2}\ h^2\,.
\eean
Consequently,
\beqn
\label{xx1}
\mathcal{L} h \le 0 \;\;\mbox{ in }\;\; \RR^N\times (0,\infty)\,,
\eeqn
where the parabolic differential operator $\mathcal{L}$ is given by 
$$
\mathcal{L} z := z_t - \Delta z + q\ \vert\nabla u\vert^{q-2}\ \nabla u \cdot
\nabla z + {q\ (q-1)\ \vert\nabla u\vert^{q-2}\over \vert\xi\vert^2}\ z^2
\,.
$$

On the one hand, since $q\in (1,2]$ and $\vert\nabla u(x,t)\vert \le
\Vert\nabla u_0\Vert_{\infty}$, it is straightforward to check that 
$$
H_1(t) := \left( {1\over \Vert h(0)\Vert_{\infty}} + {q\ (q-1)\ t\over
\vert\xi\vert^2\ \Vert\nabla u_0\Vert_{\infty}^{2-q}} \right)^{-1}\,,
\qquad t>0\,,
$$
satisfies $\mathcal{L}H_1\ge 0$ with $H_1(0)\ge h(x,0)$ for all
$x\in\RR^N$. The 
comparison principle then entails that $h(x,t)\le H_1(t)$ for $(x,t)\in
\RR^N\times (0,\infty)$, from which we conclude that 
$$
h(x,t) \le \Vert h(0)\Vert_{\infty} \le \Vert D^2 u_0\Vert_{\infty}\
\vert\xi\vert^2\,,
$$
and
$$
h(x,t) \le {\vert\xi\vert^2\ \Vert\nabla u_0\Vert_{\infty}^{2-q}\over q\
(q-1)\ t}\,.
$$
In other words, (\ref{hess1}) and (\ref{hess3}) hold true. 

On the other hand, we infer from (\ref{grad}) that 
$$
H_2(t) := {2\ C_1^{2-q}\ \vert\xi\vert^2 \over q^2\ (q-1)}\ {\Vert
u_0\Vert_{\infty}^{(2-q)/q}\over t^{2/q}}\,, \qquad t>0\,,
$$
satisfies $\mathcal{L} H_2\ge 0$ with $H_2(0)=+\infty \ge h(x,0)$ for
all $x\in\RR^N$. We then use again the comparison principle as above
and obtain (\ref{hess2}). \qed

\medskip

\begin{remark}
Since $q\in (1,2]$ and $\nabla u$ may vanish, the proof of 
Proposition~\ref{prestlap} is somehow formal because of the negative 
powers of $\vert\nabla u\vert$ in \rf{prop:2}. It can be made rigorous 
by first considering the regularised equation
$$
u_t^\varepsilon - \Delta u^\varepsilon 
+ \left( \left\vert \nabla u^\varepsilon \right\vert^2 
+ \varepsilon^2 \right)^{p/2} = 0
$$
for $\varepsilon\in (0,1)$, and then letting $\varepsilon\to 0$ as in 
\cite{BL99}.
\end{remark}

\medskip

In fact, we need a particular case of Proposition~\ref{prestlap}.

\begin{cor}\label{vittel}
Under the assumptions of Proposition~\ref{prestlap}
\bear
\label{lap1}
\Delta u(x,t) & \le & {C_3\ \Vert\nabla u_0\Vert_{\infty}^{2-q}\over
  t}\,, \\ 
\label{lap2}
\Delta u(x,t) & \le & {C_4\ \Vert u_0\Vert_{\infty}^{(2-q)/q}\over
t^{2/q}}\,,
\eear
for $(x,t)\in\RR^N\times (0,+\infty)$, where $C_3$ and $C_4$ are
positive constants depending only on $q$ and $N$. 

Furthermore, if $u_0\in W^{2,\infty}(\RR^N)$, 
\beqn
\label{lap3}
\sup_{x\in\RR^N} \Delta u(x,t) \le \sup_{x\in\RR^N} \Delta u_0(x)\,,
\quad t\ge 0\,. 
\eeqn
\end{cor}

\noindent{\sc Proof.} Consider $i\in\{1,\ldots,N\}$ and define
$\xi^i=(\xi_j^i)\in\RR^N$ by $\xi_i^i=1$ and $\xi_j^i=0$ if $j\ne
i$. We take $\xi=\xi^i$ in \rf{xx1} and obtain that $\mathcal{L}
u_{x_i x_i} \le 0$, that is, 
$$
\left( u_{x_i x_i} \right)_t - \Delta u_{x_i x_i} +q\ \vert\nabla
u\vert^{q-2}\ \nabla u . \nabla u_{x_i x_i} + q\ (q-1)\ \vert\nabla
u\vert^{q-2}\ u_{x_i x_i}^2\le 0 
$$
in $\RR^N\times (0,\infty)$. Summing the above inequality over
$i\in\{1,\ldots,N\}$ and recalling that 
$$
\vert\Delta u\vert^2 \le N\ \sum_{i=1}^N u_{x_i x_i}^2\,,
$$
we end up with 
$$
\left( \Delta u \right)_t - \Delta \left( \Delta u \right) + q\
\vert\nabla u\vert^{q-2}\ \nabla u . \nabla \left( \Delta u \right) +
{q\ (q-1)\ \vert\nabla u\vert^{q-2}\over N}\ \vert\Delta u\vert^2\le 0 
$$
in $\RR^N\times (0,\infty)$. We next proceed as in the proof of
Proposition \ref{prestlap} to complete the proof of Corollary
\ref{vittel}. 
\qed

%%%%%%%%%%%%%%%%%%%%%%%%%%%%%%%%%%%%%%%%%%%%%%%%%%%%%%%%%%%%%%%%%%%%%%%%%%
%%%%
%%%%%%%%%%%%%%%%%%%%%%%%%%%%%%%%%%%%%%%%%%%%%%%%%%%%%%%%%%%%%%%%%%%%%%%%%%
%%%%%
\section{Diffusion-dominated case}
\setcounter{equation}{0}
\setcounter{remar}{0}

The proofs of Theorems~\ref{th:01} and~\ref{th:03} rely on some
properties of the non-homogeneous heat equation which we state
now. Similar results have already been used in \cite{BGK03,LS03}.

\begin{theo} \label{tw:L1}
Assume that $u=u(x,t)$ is the solution of the Cauchy problem to the 
linear non-homogeneous heat equation 
\begin{eqnarray}  
&&u_t=\Delta u +f(x,t), \quad x\in\RR^N,\;t>0,\label{eq-n}\\
&&u(x,0)=u_0(x), \quad x\in\RR^N, \label{ini-n} 
\end{eqnarray}
with $u_0\in L^1(\RR^N)$ and $f\in L^1(\RR^N\times (0,\infty))$.
Then
\begin{equation}
\lim_{t\to\infty} \|u(t)-I_\infty G(t)\|_1 = 0\,,
\label{MG1}
\end{equation}
where
$$
I_\infty := \lim_{t\to\infty}\int_{\RR^N}u(x,t)\;dx
=\int_{\RR^N}u_0(x)\;dx +\int_0^\infty\!\! \int_{\RR^N}f(x,t)\;dx\, dt\,.
$$

Assume further that there is $p\in [1,\infty]$ such that $f(t)\in
L^p(\RR^N)$ for every $t>0$ and 
\beqn
\label{luchon}
\lim_{t\to \infty} t^{1+(N/2)(1-1/p)}\ \Vert f(t)\Vert_p=0\,. 
\eeqn
Then
\begin{equation}
\lim_{t\to \infty} t^{(N/2)(1-1/p)}\|u(t)-I_\infty G(t)\|_p = 0\,,\label{MGp}
\end{equation}
and
\begin{equation}
\lim_{t\to \infty} t^{(N/2)(1-1/p)+1/2}\|\nabla u(t)-I_\infty \nabla
G(t)\|_p = 0 \,. 
\label{MG2}
\end{equation}
\end{theo}

\proof
We first observe that the assumptions on $u_0$ and $f$ warrant that
$I_\infty$ is finite, and we refer to \cite{BGK03} for the proof of
\rf{MG1}. We next assume \rf{luchon} and prove \rf{MG2}. Let $T>0$ and
$t\in (T,\infty)$. By the Duhamel formula,  
$$
\nabla u(t)=\nabla G(t-T)*u(T) +\int_T^t  \nabla G(t-\tau)*f(\tau)\;d\tau.
$$  
It follows from the Young inequality that
\bean
& & t^{(N/2)(1-1/p)+1/2}\ \Vert \nabla u(t) - \nabla G(t-T)*u(T)\Vert_p \\ 
& \le & C\ t^{(N/2)(1-1/p)+1/2}\ \int_T^{(T+t)/2}
(t-\tau)^{-(N/2)(1-1/p)-1/2}\ \Vert f(\tau)\Vert_1\ d\tau \\ 
& + & C\ t^{(N/2)(1-1/p)+1/2}\ \int_{(T+t)/2}^t (t-\tau)^{-1/2}\ \Vert
f(\tau)\Vert_p\ d\tau \\  
& \le & C\ \left( \frac{t}{t-T} \right)^{(N/2)(1-1/p)+1/2}\
\int_T^\infty \Vert f(\tau)\Vert_1\ d\tau \\ 
& + & C\ \sup_{\tau\ge T} \left\{ \tau^{(N/2)(1-1/p)+1}\ \Vert
  f(\tau)\Vert_p \right\}\ \int_{(T+t)/2}^t (t-\tau)^{-1/2}\
\tau^{-1/2}\ d\tau \\ 
& \le & C\ \left( \frac{t}{t-T} \right)^{(N/2)(1-1/p)+1/2}\
\int_T^\infty \Vert f(\tau)\Vert_1\ d\tau \\ 
& + & C\ \sup_{\tau\ge T} \left\{ \tau^{(N/2)(1-1/p)+1}\ \Vert
  f(\tau)\Vert_p \right\}\,. 
\eean
Also, classical properties of the heat semigroup (see, e.g.,
\cite{DZ92}) ensure that  
$$
 \lim_{t\to\infty} t^{(N/2)(1-1/p)+1/2} \left\|\nabla
 G(t-T)*u(T)-\left(\int_{\RR^N} u(x,T)\;dx\right) \nabla
 G(t-T)\right\|_p = 0\,, 
$$
and
$$ 
 \lim_{t\to\infty} t^{(N/2)(1-1/p)+1/2} \left\|\nabla G(t-T)-\nabla
 G(t)\right\|_p = 0 
$$
for every $p\in [1,\infty]$.
Since, by elementary calculations, we have  
\bean
& & \Vert \nabla u(t) - I_\infty\ \nabla G(t)\Vert_p \\
& \le & \Vert \nabla u(t) - \nabla G(t-T)*u(T)\Vert_p \\
& + & \left\Vert \nabla G(t-T)*u(T) - \left( \int_{\RR^N} u(x,T)\, dx
  \right)\ \nabla G(t-T) \right\Vert_p \\ 
& + & \left\vert \int_{\RR^N} u(x,T)\, dx - I_\infty \right\vert\
\Vert\nabla G(t-T)\Vert_p + \vert I_\infty\vert\ \Vert \nabla
G(t-T)-\nabla G(t)\Vert_p\,, 
\eean
the previous relations imply that
\bean
& & \limsup_{t\to\infty}\, t^{(N/2)(1-1/p)+1/2}\ \Vert \nabla u(t) -
I_\infty\ \nabla G(t)\Vert_p \\ 
& \le & C\ \left( \int_T^\infty \Vert f(\tau)\Vert_1\ d\tau +
  \sup_{\tau\ge T} \left\{ \tau^{(N/2)(1-1/p)+1}\ \Vert f(\tau)\Vert_p
  \right\} + \left\vert \int_{\RR^N} u(x,T)\, dx - I_\infty
  \right\vert \right)\,. 
\eean
The above inequality being valid for any $T>0$, we may let $T\to
\infty$ and conclude that (\ref{MG2}) holds true. The assertion
(\ref{MGp}) then follows from (\ref{MG1}) and (\ref{MG2}) by the
Gagliardo-Nirenberg inequality. \qed 

\bigskip

\noindent
{\sc Proof of Theorem \ref{th:01}.} Since $u$ is non-negative, we
infer from \cite[Eq.~(17)]{BL99} that there is a constant $C=C(q)$
such that  
$$
\Vert \nabla u^{(q-1)/q}(t) \Vert_\infty \le C\ \Vert
u(t/2)\Vert_\infty^{(q-1)/q}\ t^{-1/2}\,, \quad t>0\,. 
$$
Also, $u$ is a subsolution to the linear heat equation and therefore satisfies
$$
\Vert u(t)\Vert_p \le \Vert G(t)*u_0\Vert_p \le C\ t^{-(N/2)(1-1/p)}\
\Vert u_0\Vert_1 \,, \quad t>0\,,
$$
for every $p\in [1,\infty]$ by the comparison principle. Since $\nabla
u = (q/(q-1))\ u^{1/q}\ \nabla u^{(q-1)/q}$, we obtain that  
$$
t^{(N/2)(1-1/p)+1}\ \Vert \vert \nabla u(t)\vert^q \Vert_p \le C\
t^{(N+2-q(N+1))/2} \mathop{\longrightarrow}_{t\to\infty} 0 
$$
for $p\in [1,\infty]$, because $q>(N+2)/(N+1)$. Theorem \ref{th:01}
then readily follows by Theorem \ref{tw:L1} with $f(x,t)=-|\nabla
u(x,t)|^q$. \qed 

\bigskip

\noindent
{\sc Proof of Theorem \ref{th:03}, part a).}
Since $q\ge 2$, we infer from \cite{LS03} that $I_\infty$ is finite
and negative and that  
\beqn
\label{t1rex}
\nabla u \in L^q(\RR^N\times (0,\infty))\,.
\eeqn
Setting $b:= \Vert\nabla u_0\Vert_\infty^{q-2}$, it follows from
\rf{contrex} that $u_t-\Delta u \ge -b\ \vert\nabla u\vert^2$ in
$\RR^N\times (0,\infty)$. The comparison principle then entails that
$u\ge w$, where $w$ is the solution to  
$$
w_t - \Delta w = -b\ \vert\nabla w\vert^2\,, \quad w(0)=u_0\,.
$$
The Hopf-Cole transformation $h:=e^{-bw}-1$ then implies that $h$ solves 
$$
h_t - \Delta h = 0\,, \quad h(0)=e^{-bu_0}-1\,.
$$
Therefore, for $t>0$, 
$$0\le -b w(x,t) \le h(x,t) \le \Vert h(t)\Vert_\infty \le C\
t^{-N/2}\ \Vert h(0)\Vert_1 \le C\ t^{-N/2}\,,$$  
since $u_0\in L^1(\RR^N)\cap L^\infty(\RR^N)$. Recalling that $0\ge u
\ge w$, we end up with  
\beqn
\label{t2rex}
\Vert u(t)\Vert_\infty \le C\ t^{-N/2}\,, \quad t>0\,.
\eeqn
It next follows from \cite[Theorem~2]{GGK03} that 
$$
\Vert\nabla u(t)\Vert_\infty \le C\ \Vert u(t/2)\Vert_\infty \
t^{-1/2}\,, \quad t>0\,, 
$$
which, together with (\ref{t2rex}), yields
\beqn
\label{t3rex}
\Vert\nabla u(t)\Vert_\infty \le C\ t^{-(N+1)/2}\,, \quad t>0\,.
\eeqn
Recalling \rf{contrex}, we also have
\beqn
\label{adrenaline}
\Vert\nabla u(t)\Vert_\infty \le C\ (1+t)^{-(N+1)/2}\,, \quad t\ge
0\,.
\eeqn

We next put 
$$
\mathcal{A}_1(t) := \sup_{\tau\in (0,t)} \left\{ \tau^{1/2}\
  \Vert\nabla u(\tau)\Vert_1 \right\}\,, 
$$
which is finite by \cite{BSW02}. Since $q\ge 2$ and $N\ge 1$, we infer
from the Duhamel formula and (\ref{adrenaline}) that, for $\alpha\in
(0,1/2)$,   
\bean
t^{1/2}\ \Vert\nabla u(t)\Vert_1 & \le & C\ \Vert u_0\Vert_1 + C\
t^{1/2}\ \int_0^t (t-\tau)^{-1/2}\ \Vert\nabla u(\tau)\Vert_q^q\ d\tau \\ 
& \le & C + C\ t^{1/2}\ \int_0^t (t-\tau)^{-1/2}\
(1+\tau)^{-(q-1)(N+1)/2}\ \Vert\nabla u(\tau)\Vert_1\ d\tau \\
& \le & C + C\ t^{1/2}\ \int_0^t (t-\tau)^{-1/2}\
(1+\tau)^{-1}\ \tau^{-1/2}\ \mathcal{A}_1(\tau)\ d\tau \\ 
& \le & C + C\ \alpha^{-1/2}\ \int_0^{(1-\alpha) t} (1+\tau)^{-1}\
\tau^{-1/2} \mathcal{A}_1(\tau)\ d\tau \\  
& + & C\ t^{1/2}\ \mathcal{A}_1(t)\ \int_{(1-\alpha) t}^t
(t-\tau)^{-1/2}\ \frac{2}{2+t}\ \tau^{-1/2}\ d\tau\\ 
& \le & C + C\ \alpha^{-1/2}\ \int_0^t
(1+\tau)^{-1}\ \tau^{-1/2} \mathcal{A}_1(\tau)\ d\tau \\ 
& + & C\ \mathcal{A}_1(t)\ \int_{1-\alpha}^1 (1-\tau)^{-1/2}\
\tau^{-1/2}\ d\tau\,,
\eean
whence
$$
\left( 1 - C\ \alpha^{1/2} \right)\ \mathcal{A}_1(t) \le C(\alpha)\ \left(
  1+\int_0^t (1+\tau)^{-1}\ \tau^{-1/2} \mathcal{A}_1(\tau)\
  d\tau \right)\,. 
$$
Consequently, there is $\alpha_0\in (0,1/2)$ sufficiently small such that 
$$
\mathcal{A}_1(t) \le \mathcal{B}_1(t) := C(\alpha_0)\ \left( 1+\int_0^t
  (1+\tau)^{-1}\ \tau^{-1/2} \mathcal{A}_1(\tau)\ d\tau \right) 
$$
for $t\ge 0$. Now, for $t\ge 0$, 
$$
{d\mathcal{B}_1\over dt}(t) = C(\alpha_0)\ (1+t)^{-1}\ t^{-1/2}
\mathcal{A}_1(t) \le C(\alpha_0)\ (1+t)^{-1}\ t^{-1/2}
\mathcal{B}_1(t)\,,  
$$
from which we deduce that
$$
\mathcal{A}_1(t) \le \mathcal{B}_1(t) \le \mathcal{B}_1(0)\
\exp{\left\{ C(\alpha_0)\ \int_{0}^t (1+\tau)^{-1}\
    \tau^{-1/2}\ d\tau \right\}} \le C(\alpha_0)\,.
$$
We have thus proved that
\beqn
\label{t4rex}
\Vert \nabla u(t)\Vert_1\le C\ t^{-1/2}\,, \quad t>0\,.
\eeqn

We finally infer from (\ref{t3rex}), (\ref{t4rex}) and the H\"older
inequality that  
$$
t^{(N/2)(1-1/p)+1}\ \Vert \vert \nabla u(t)\vert^q \Vert_p \le C\
t^{(N+2-q(N+1))/2} \mathop{\longrightarrow}_{t\to\infty} 0 
$$
for $p\in [1,\infty]$, and we conclude as in the proof of Theorem
\ref{th:01}. \qed 

\bigskip

\noindent
{\sc Proof of Theorem \ref{th:03}, part b).}
Since $q\in (q_c,2)$, we obtain from \cite{LS03} that there is
$\varepsilon>0$ such that, if $u_0$ fulfils \rf{small:u0}, then
$I_\infty$ is finite and negative and there are $C>0$ and $\delta>0$
such that   
\beqn
\label{t5rex}
\Vert \nabla u(t)\Vert_q^q \le C\ t^{-1}\ (1+t)^{-\delta}\,, \quad t>0\,.
\eeqn
In particular, 
\beqn
\label{t6rex}
\vert\nabla u\vert^q \in L^1(\RR^N\times (0,\infty)) \;\mbox{ and
}\;\lim_{t\to\infty} t\ \left\Vert \vert\nabla u(t)\vert^q
\right\Vert_1 = 0\,. 
\eeqn

We next claim that 
\beqn
\label{t7rex}
\Vert \nabla u(t) \Vert_\infty \le C\ t^{-(N+1)/2}\,, \quad t>0\,.
\eeqn
Indeed, we fix $r\in (q_c,q)$ such that $r<N/(N-1)$ and define
$s=r/(r-1)$ and a sequence $(r_i)_{i\ge 0}$ by $$ 
r_0={1\over q} \;\mbox{ and }\; r_{i+1} = {(N+1)\ r - (N+2) \over 2r}
+ {q\over r}\ r_i\,, \quad i\ge 0\,. 
$$
We now proceed by induction to show that, for each $i\ge 0$, there is
$K_i\ge 0$ such that
\beqn
\label{t8rex}
\Vert\nabla u(t)\Vert_\infty \le K_i\ \left( t^{-(N+1)/2} + t^{-r_i}
\right)\,, \quad t>0\,.
\eeqn
Thanks to (\ref{grad}), the assertion \rf{t8rex} is true for
$i=0$. Assume next that \rf{t8rex} holds true for some $i\ge 0$. We
infer from \rf{t5rex}, \rf{t8rex} and the Duhamel formula that
\bean
\Vert\nabla u(t)\Vert_\infty & \le & C\ \Vert u_0\Vert_1\ t^{-(N+1)/2}
+ C\ \int_0^{t/2} (t-\tau)^{-(N+1)/2}\ \Vert\nabla u(\tau)\Vert_q^q\
d\tau \\
& + & C\ \int_{t/2}^t (t-\tau)^{-(N/2)(1-1/r)-1/2}\ \Vert\nabla
u(\tau)\Vert_{sq}^q\ d\tau \\
& \le & C\ t^{-(N+1)/2}\ \left( \Vert u_0\Vert_1 + \int_0^{t/2}
  \Vert\nabla u(\tau)\Vert_q^q\ d\tau \right) \\
& + &  C\ \int_{t/2}^t (t-\tau)^{-(N/2)(1-1/r)-1/2}\ \Vert\nabla
u(\tau)\Vert_\infty^{q/r}\ \Vert\nabla u(\tau)\Vert_{q}^{q/s}\ d\tau \\
& \le & C\ t^{-(N+1)/2} + C\ \mathcal{I}(t)\,,
\eean
where
$$
\mathcal{I}(t) := \int_{t/2}^t (t-\tau)^{-(N/2)(1-1/r)-1/2}\ \left(
  \tau^{-(N+1)/2} + \tau^{-r_i} \right)^{q/r}\ \tau^{-1/s}\ d\tau \,.
$$
Since $r<N/(N-1)$ and $q>q_c$, we have
\bean
\mathcal{I}(t) & \le & C\ \int_{t/2}^t (t-\tau)^{-(N/2)(1-1/r)-1/2}\ \left(
  \tau^{-(q(N+1))/2r} + \tau^{-(q r_i)/r} \right) \tau^{-1/s}\ d\tau \\
  & \le & C\ t^{-((N+1)r-(N+2))/2r}\ \left( t^{-(q(N+1))/2r} + t^{-(q
      r_i)/r}\right) \\ 
& \le & C\ \left( t^{-(N+1)/2}\ t^{-((N+1)q-(N+2))/2r} + t^{-r_{i+1}}
\right) \\ 
& \le & C\ \left( t^{-(N+1)/2} + t^{-r_{i+1}} \right)
\eean
for $t\ge 1$. Consequently, for $t\ge 1$, 
$$
\Vert\nabla u(t)\Vert_\infty \le K_{i+1}\ \left( t^{-(N+1)/2} +
  t^{-r_{i+1}} \right)\,, 
$$
while (\ref{contrex}) implies that the same inequality is valid for
$t\in [0,1]$ for  
a possibly larger constant $K_{i+1}$. Thus \rf{t8rex} is true for
$i+1$, which completes  
the proof of (\ref{t8rex}). To obtain \rf{t7rex}, it suffices to note
that $r_i\to\infty$  
since $q>r$.  

Now, owing to \rf{t6rex} and \rf{t7rex}, we are in a position to apply
Theorem \ref{tw:L1}  
and conclude that \rf{lim:G1} and \rf{lim:G2} holds true for $p=1$ and
$p=\infty$. The general case $p\in (1,\infty)$ then follows by the
H\"older inequality. \qed

%%%%%%%%%%%%%%%%%%%%%%%%%%%%%%%%%%%%%%%%%%%%%%%%%%%%%%%%%%%%%%%%%%%%%%%%%%
%%%%%
%%%%%%%%%%%%%%%%%%%%%%%%%%%%%%%%%%%%%%%%%%%%%%%%%%%%%%%%%%%%%%%%%%%%%%%%%%
%%%%%

\section{Convergence towards very singular solutions}
\setcounter{equation}{0}
\setcounter{remar}{0}

The goal of this section is to prove Theorem \ref{th:02}.
Recall that we assume that $1<q<q_c$ and that $u_0$ is a
non-negative and integrable function satisfying in addition
\beqn
\label{gru0}
\mathop{\mbox{ess lim}}_{\vert x\vert \to \infty} \vert x\vert^a\ u_0(x)
= 0\,,
\eeqn
with $a=(2-q)/(q-1)\in (N,\infty)$. We define 
$$
R(u_0) := \inf{\left\{ R>0\,, \;\; \vert x\vert^a\ u_0(x) \le \gamma_q
\;\mbox{ a.e. in}\; \{\vert x\vert\ge R\} \right\}}\,,
$$
where $\gamma_q := (q-1)^{(q-2)/(q-1)}\ (2-q)^{-1}$ and observe that
$R(u_0)$ is finite by (\ref{gru0}). 

Denoting by $u$ the corresponding solution to (\ref{eq}) and introducing 
$$
\tau(u_0) := \left( {(N+2)-q (N+1)\over (N+1) q - N} \right)^{1-q}\
R(u_0)^2\,,
$$
we infer from \cite[Lemma 2.2 \& Proposition 2.4]{BL01} that there is a
constant $C_1$ depending only on $N$ and $q$ such that
\beqn
\label{bdu}
t^{(a-N)/2}\ \Vert u(t)\Vert_{1} + t^{a/2}\ \Vert u(t)\Vert_{\infty} +
t^{(a+1)/2}\ \Vert \nabla u(t)\Vert_{\infty} \le C_1
\eeqn
for each $t>\tau(u_0)$ and
\beqn
\label{bsup}
u(x,t) \le \Gamma_q(\vert x\vert - R(u_0))\,, \qquad t>0\,, \;\;\; \vert
x\vert>R(u_0)\,.
\eeqn
Here, $\Gamma_q$ is given by $\Gamma_q(r) = \gamma_q\ r^{-a}$, $r\in
(0,\infty)$.

\medskip

Let us observe at this point that decay estimates for $\nabla u(t)$ in
$L^p$ can be deduced from (\ref{bdu}) and the Duhamel formula. 

\begin{lem}\label{wattwiller}
For $p\in [1,\infty]$, there is a constant $C(p)$ depending only on $N$, $q$ 
and $p$ such that
\beqn
\label{estgup}
t^{((a+1)p-N)/2p}\ \Vert\nabla u(t)\Vert_{L^p} \le C(p) \;\;\mbox{ for
}\;\; t>\tau(u_0)\,. 
\eeqn
\end{lem}

\proof
Indeed, since $u$ is non-negative, it follows from
\cite[Theorem~1]{BL99} that 
$$
\Vert\nabla u^{(q-1)/q}(t)\Vert_\infty \le C(q)\ t^{-1/q} 
$$
for $t>0$, which, together with (\ref{bdu}) and the
Duhamel formula  entails that, for $t>\tau(u_0)$, 
\bean
\Vert\nabla u(t)\Vert_{1} & \le & \Vert\nabla G(t/2) * u(t/2)\Vert_{1}
+ \int_{t/2}^t \left\Vert \nabla G(t-s) * \vert\nabla u\vert^q \right\Vert_{1}\
ds \\
& \le & C\ t^{-1/2}\ \Vert u(t/2)\Vert_{1} + C\ \int_{t/2}^t
(t-s)^{-1/2}\ \left\Vert \nabla u^{(q-1)/q}(s) \right\Vert_{\infty}^q\
\Vert u(s)\Vert_{1}\ ds \\
& \le & C\ t^{-(a+1-N)/2} + C\ \int_{t/2}^t (t-s)^{-1/2}\ s^{-(a+2-N)/2}\
ds \\
& \le & C\ t^{-(a+1-N)/2}\,.
\eean
Interpolating between (\ref{bdu}) and the above estimate yields
\rf{estgup}. \qed 

\bigskip

In order to investigate the large time behavior of $u$, we use a
rescaling method and introduce the sequence of rescaled solutions
$(u_k)_{k\ge 1}$ defined by
$$
u_k(x,t) = k^a\ u(k x, k^2 t)\,, \qquad (x,t)\in \RR^N\times
[0,\infty)\,, \;\;\; k\ge 1\,.
$$

\begin{lem}\label{levss1} 
For $k\ge 1$, we have
\beqn
\label{bdku}
t^{(a-N)/2}\ \Vert u_k(t)\Vert_{1} + t^{a/2}\ \Vert
u_k(t)\Vert_{\infty} + t^{(a+1)/2}\ \Vert \nabla u_k(t)\Vert_{\infty}
\le C_1
\eeqn
for $t>\tau_k := \tau(u_0)\ k^{-2}$ and
\beqn
\label{bksup}
u_k(x,t) \le \Gamma_q\left( \vert x\vert - {R(u_0)\over k} \right)
\;\;\mbox{ for }\;\; \vert x\vert>{R(u_0)\over k} \;\mbox{ and }\; t>0\,.
\eeqn
\end{lem}

\proof
It is  straightforward to check that, for each $k\ge 1$, $u_k$ is the
solution to (\ref{eq}) with initial datum $u_k(0)$ and satisfies  estimates 
\rf{bdku} and \rf{bksup} as a consequence of (\ref{bdu}) and (\ref{bsup}). 
\qed

\bigskip

We next use (\ref{eq}) and the non-negativity of $u_k$ to control the
behavior of $u_k(x,t)$ for large $x$ uniformly with respect to $k$. For 
$k\ge 1$, $t>0$ and $R\ge 0$, we put
\beqn
\label{aix}
I_k(R,t) := \int_{\{ \vert x \vert \ge R\}} u_k(x,t)\ dx + \int_0^t \int_{\{
\vert x \vert \ge R\}} \vert\nabla u_k(x,t)\vert^q\ dxdt\,.
\eeqn

\begin{lem}\label{levss2}
For every $T>0$, we have
\beqn
\label{uik}
\lim_{R\to \infty} \sup_{k\ge 1}\, \sup_{t\in [0,T]}\, I_k(R,t) = 0\,. 
\eeqn
\end{lem}

\proof
 Let $\varrho$ be a non-negative function in $\mathcal{C}^\infty(\RR^N)$
 such that $0\le \varrho \le 1$ and 
$$
\varrho(x) = 0 \;\;\mbox{ if }\;\; \vert x\vert \le {1\over 2} \;\;\mbox{
and }\;\; \varrho(x) = 1 \;\;\mbox{ if }\;\; \vert x\vert \ge 1\,.
$$
For $R>0$ and $x\in\RR^N$, we set $\varrho_R(x) = \varrho(x/R)$. As $u_k$
is a non-negative solution to (\ref{eq}), we have
\bear
I_k(R,t) & \le & \int u_k(x,t)\ \varrho_R(x)\ dx + \int_0^t \int \vert\nabla
u_k(x,s)\vert^q\ \varrho_R(x)\ dxds \nonumber\\
& \le & \int u_k(x,0)\ \varrho_R(x)\ dx + \int_0^t \int u_k(x,s)\
\vert\Delta\varrho_R(x)\vert\ dxds \nonumber\\
& \le & k^{a-N}\ \int_{\{ \vert x\vert \ge kR/2\}} u_0(x)\ dx +
{\vert\Delta\varrho\vert_{\infty}\over R^2}\ \int_0^t \int_{\{R/2\le
\vert x\vert\le R\}} u_k(x,s)\ dxds\,. \label{xx3}
\eear
Owing to (\ref{gru0}) and (\ref{bksup}), we further obtain that, for $R\ge
1+4\ R(u_0)$, 
\bean
I_k(R,t) & \le & k^{a-N}\ \int_{\{ \vert x\vert \ge kR/2\}} \Gamma_q\left(
{\vert x\vert\over 2} \right)\ dx \\
&& +\ 
{\vert\Delta\varrho\vert_{\infty}\over R^2}\ \int_0^t \int_{\{R/2\le
\vert x\vert\le R\}} \Gamma_q\left( \vert x\vert - {R(u_0)\over k}
\right)\ dxds \\
& \le & C\ R^{-(a-N)} + {T\ \vert\Delta\varrho\vert_{\infty}\over R^2}\
\int_{\{R/2\le \vert x\vert\le R\}} \Gamma_q\left( {R\over 4} \right)\ dx
\\
& \le & C(T,\varrho)\ R^{-(a-N)}\,. 
\eean
Lemma~\ref{levss2} then readily follows since $a>N$. \hfill $\square$

\medskip

We finally study the behavior of $u_k$ for small times. 

\begin{lem}\label{levss3} 
Let $r>0$. There is a positive constant $C(r)$ depending only on $q$, $N$
and $r$ such that
\beqn
\label{smallt}
\int_{\{\vert x\vert \ge r\}} u_k(x,t)\ dx \le C(r)\ \left( \sup_{\vert
x\vert\ge k r/2} \left\{ \vert x\vert^a\ u_0(x) \right\} + t \right)
\eeqn
for $t>\tau_k$ and $k\ge 4\ R(u_0)/r$.
\end{lem}

\proof
 We fix $r>0$ and use the same notations as in the proof of Lemma
 \ref{levss2}. Thanks to the properties of $\varrho$, we infer from
 (\ref{xx3}) with $R=r$ that, for $t>\tau_k$ and $k\ge 4\ R(u_0)/r$, 
\bean
\int_{\{\vert x\vert \ge r\}} u_k(x,t)\ dx & \le & \int u_k(x,t)\
\varrho_r(x)\ dx \\
& \le & k^{a-N}\ \int_{\{ \vert x\vert \ge k r/2\}} u_0(x)\ dx \\
& + & {\vert\Delta\varrho\vert_{\infty}\over r^2}\ \int_0^t 
\int_{\{r/2\le \vert x\vert\le r\}} u_k(x,s)\ dxds \\
& \le & C(\varrho,r)\ \left( \sup_{\vert x\vert\ge k r/2} \left\{ \vert
x\vert^a\ u_0(x) \right\} + t \right)\,,
\eean
where we have used (\ref{bksup}) to obtain the last inequality. \hfill
$\square$

\medskip

\noindent{\sc Proof of Theorem~\ref{th:02}.} Owing to
Lemma~\ref{levss1} and Lemma~\ref{levss2} we may proceed as in
\cite[Theorem 3]{BL99} to prove that 
there are a subsequence of $(u_k)$ (not relabeled) and a non-negative
function 
$$
u_\infty \in \mathcal{C}((0,\infty);L^1(\RR^N)) \cap
L^q((s,\infty)\times\RR^N)) \cap L^\infty(s,\infty;W^{1,\infty}(\RR^N))
$$
satisfying
$$
u_\infty(t) = G(t-s) * u_\infty(s) - \int_s^t G(t-\tau) * \vert\nabla
u_\infty(\tau)\vert^q\ d\tau
$$
and
\beqn
\label{cvk}
\lim_{k\to \infty}\ \sup_{\tau\in [s,t]} \Vert u_k(\tau) -
u_\infty(\tau)\Vert_{1} = 0
\eeqn
for every $s>0$ and $t>s$.

It remains to identify the behavior of $u_\infty$ as $t\to 0$. On the one
hand, consider $r>0$ and $t>0$. Since $\tau_k\to 0$ as $k\to \infty$, we
have $t>\tau_k$ for $k$ large enough and it follows from Lemma
\ref{levss3}, (\ref{gru0}) and (\ref{cvk}) that 
$$
0\le \int_{\{\vert x\vert \ge r\}} u_\infty(x,t)\ dx \le C(r)\ t\,.
$$
Consequently,
\beqn
\label{xx4}
\lim_{t\to 0} \int_{\{\vert x\vert \ge r\}} u_\infty(x,t)\ dx = 0\,. 
\eeqn

On the other hand, consider $M>0$ and set $k_M:=M^{1/(a-N)}$. For $k\ge
k_M$, we denote by $v_k$ the solution to (\ref{eq}) with initial datum
$v_k(0)$ given by $v_k(x,0):= M\ k^N\ u_0(k x)$, $x\in\RR^N$. Since $a>N$,
we have $v_k(0)\le u_k(0)$ for $k\ge k_M$ and the comparison principle
warrants that
\beqn
\label{comp}
v_k(x,t)\le u_k(x,t)\,, \qquad (x,t)\in \RR^N\times [0,\infty)\,, \qquad
k\ge k_M\,.
\eeqn
We next observe that $(v_k(0))$ converges narrowly towards $(M\ \Vert
u_0\Vert_{1})\ \delta$ as $k\to \infty$ ($\delta$ denoting the Dirac
mass at $x=0$). We then proceed as in \cite{BL99} to conclude that 
$$
\lim_{k\to \infty}\ \sup_{\tau\in [s,t]} \Vert v_k(\tau) -
S_{M}(\tau)\Vert_{1} = 0
$$
for every $s>0$ and $t>s$, where $S_{M}$ denotes the unique non-negative
solution to (\ref{eq}) with initial datum $(M\ \Vert u_0\Vert_{1})\
\delta$ \cite{BL99}. Recalling (\ref{cvk}) and (\ref{comp}), we realize
that 
$$
S_{M}(x,t) \le u_\infty(x,t)\,, \qquad (x,t)\in \RR^N\times (0,\infty)\,. 
$$
The above inequality being valid for any $M>0$, it is then straightforward
to deduce that
\beqn
\label{xx5}
\lim_{t\to 0} \int_{\{\vert x\vert \le r\}} u_\infty(x,t)\ dx = 
\infty\,. 
\eeqn

In other words, $u_\infty$ is a very singular solution to (\ref{eq}) and
the uniqueness of the very singular solution to (\ref{eq})
(cf.~\cite{BKLxx,QW01}) implies that $u_\infty=W$, where $W$ is the
very singular solution to \rf{eq}, see Theorem~\ref{th:02}. The
uniqueness of the limit 
actually entails that the whole sequence $(u_k)_{k\ge 1}$ converges
towards $W$ in $\mathcal{C}([s,t];L^1(\RR^N))$ for $s>0$ and $t>s$.
Expressed in terms of $u$, we have thus shown that
\beqn
\label{cvl1}
\lim_{t\to \infty} t^{(a-N)/2}\ \Vert u(t) - W(t) \Vert_{1} = 0\,.
\eeqn
Finally, it follows from (\ref{bdu}), (\ref{cvl1}) and the
Gagliardo-Nirenberg inequality that \rf{lim:W1} holds true. 

The last step of the proof is to obtain the convergence \rf{lim:W2} for the
gradients. Consider $p\in [1,\infty]$, $t>0$ and $\alpha\in (0,1)$. By the
Duhamel formula, we have
\bean
A_p(t) & := & t^{((a+1)p-N)/2p}\ \Vert \nabla (u-W)(t)\Vert_{L^p} \\
& \le & t^{((a+1)p-N)/2p}\ \Vert \nabla G((1-\alpha) t) * (u-W)(\alpha t)
\Vert_{L^p} \\
& + & t^{((a+1)p-N)/2p}\ \int_{\alpha t}^t \left\Vert \nabla G(t-s) * 
\left( \vert\nabla u(s)\vert^q - \vert\nabla W(s)\vert^q \right)
\right\Vert_{L^p}\ ds \\
& \le & C(\alpha)\ t^{(a-N)/2}\ \Vert (u-W)(\alpha t) \Vert_{1} \\
& + & C\ t^{((a+1)p-N)/2p}\ \int_{\alpha t}^t (t-s)^{-1/2}\ s^{-1/2}\
\Vert \nabla (u-W)(s) \Vert_{L^p}\ ds\,,
\eean
where we have used the fact that
$$
\max{\left\{ \Vert\nabla u(s)\Vert_{\infty} , \Vert\nabla
W(s)\Vert_{\infty} \right\}} \le C\ s^{-(a+1)/2}
$$
by (\ref{bdu}) and the properties of $W$ in order to obtain the last
inequality.  
Consequently, by the definition of $A_p(t)$ and the change of
variables $s\mapsto ts$, we obtain 
\bean
A_p(t) & \le & C(\alpha)\ t^{(a-N)/2}\ \Vert (u-W)(\alpha t) \Vert_{1}
\\
& + & C\ t^{((a+1)p-N)/2p}\ \int_{\alpha t}^t (t-s)^{-1/2}\ s^{-1/2}\
s^{-((a+1)p-N)/2p}\ A_p(s)\ ds \\
& \le & C(\alpha)\ t^{(a-N)/2}\ \Vert (u-W)(\alpha t) \Vert_{1} \\
& + & C\ \int_\alpha^1 (1-s)^{-1/2}\ s^{-1/2}\ s^{-((a+1)p-N)/2p}\ A_p(s
t)\ ds\,.
\eean
Now, introducing 
$$
A_p(\infty) := \limsup_{t\to +\infty} A_p(t) \ge 0\,,
$$
which is finite by (\ref{estgup}), we may let $t\to +\infty$ in the above
inequality and use (\ref{cvl1}) to conclude that
$$
A_p(\infty) \le C\ \int_\alpha^1 (1-s)^{-1/2}\ s^{-1/2}\
s^{-((a+1)p-N)/2p}\ ds\ A_p(\infty)\,.
$$
Finally, the choice of $\alpha<1$ sufficiently close to 1 
 readily yields that $A_p(\infty)=0$, from which \rf{lim:W2}
follows. \qed

%%%%%%%%%%%%%%%%%%%%%%%%%%%%%%%%%%%%%%%%%%%%%%%%%%%%%%%%%%%%%%%%%%%%%%%%%%
%%%%%
%%%%%%%%%%%%%%%%%%%%%%%%%%%%%%%%%%%%%%%%%%%%%%%%%%%%%%%%%%%%%%%%%%%%%%%%%%
%%%%%

\section{Proofs of Theorems \ref{th:04} and \ref{th:05}}
\setcounter{equation}{0}
\setcounter{remar}{0}

{\sc Proof of Theorem \ref{th:04}, part a).}
The required non-positive self-similar solution 
$$
V=V(x,t)= t^{-(2-q)/(2(q-1))}V\left( x\ t^{-1/2},1 \right)
$$
is constructed and studied in \cite[Theorem~3.5]{BSW02}. In
particular, it is shown that the self-similar profile
$\mathcal{V}(x):=V(x,1)$ is a radially symmetric bounded
$\mathcal{C}^2$ function. Moreover, the profile $\mathcal{V}$ and its
first derivative $\mathcal{V}'$ both decay exponentially as
$|x|\to\infty$ (see \cite[Proposition~3.14]{BSW02}) \qed 

\bigskip

\noindent
{\sc Proof of Theorem \ref{th:04}, part b).} 
Recall that by assumption \rf{u0-}, $u=u(x,t)$ is a non-positive
solution to \rf{eq}.  
For $t\ge 0$, we put 
$m(t) = \inf{\{ u(x,t)\,, \ x\in\RR^N\}}\le 0$. The
comparison principle ensures that $t\mapsto m(t)$ is a non-decreasing
function of time and  
$$
m_\infty := \sup_{t\ge 0}\, m(t) \in (-\infty,0]\,.
$$
Since $u$ is a classical solution to (\ref{eq}), it follows from
(\ref{eq}) that 
$$
u(x,t) \le u_0(x) + \int_0^t \Delta u(x,\tau)\ d\tau \le u_0(x) + \int_0^t
\sup_{y\in\RR^N} \Delta u(y,\tau)\ d\tau 
$$
for every $x\in\RR^N$ and $t \ge 0$. Therefore, 
$$
m(t) \le - \Vert u_0\Vert_{\infty} + \int_0^t \sup_{y\in\RR^N} \Delta
u(y,\tau)\ d\tau\,,
$$
and we infer from (\ref{lap2}) and (\ref{lap3}) that
$$
m(t) \le - \Vert u_0\Vert_{\infty} + T\ \left\Vert (\Delta u_0)^+
\right\Vert_{\infty} + 
C\ \Vert u_0\Vert_{\infty}^{(2-q)/q}\ \int_T^t \tau^{-2/q}\ d\tau
$$
for $T>0$ and $t>T$. Since $q< 2$, we may let $t\to \infty$ in the
above inequality and obtain with the choice $T= \Vert
u_0\Vert_{\infty}^{(2-q)/2}\ \left\Vert (\Delta u_0)^+
\right\Vert_{\infty}^{-q/2}$ that there is a constant $K$ depending
only on $q$ such that 
\beqn
\label{xx2}
m_\infty \le - \Vert u_0\Vert_{\infty} + K^{q/2}\ \left\Vert (\Delta
  u_0)^+ \right\Vert_{\infty}^{(2-q)/2}\ \Vert
u_0\Vert_{\infty}^{(2-q)/2}\,. 
\eeqn
Therefore, if $\Vert u_0\Vert_{\infty} > K\ \left\Vert (\Delta u_0)^+
  \right\Vert_{\infty}^{(2-q)/q}$, we readily conclude from
  (\ref{xx2}) that $m_\infty< 0$, whence \rf{lim:u:infty}. \qed 

\bigskip

\noindent{\sc Proof of Theorem \ref{th:05}.} The proof of the first assertion 
of Theorem \ref{th:05} is the same as that of Theorem~\ref{th:04},
part b), hence we skip it. We next assume that $N\le 3$ and that
$1<q<4/(1+\sqrt{1+2N})$. For $t>0$, we put  
$$
\ell(t) := \Vert u(t)\Vert_\infty\ \left\Vert (\Delta u(t))^+
\right\Vert_{\infty}^{1-2/q}\,. 
$$
Since $u$ is a non-positive subsolution to the linear heat equation,
we infer from classical properties of the heat semigroup that 
$$
\Vert u(t)\Vert_\infty \ge \Vert G(t)*u_0\Vert_\infty \ge C\ t^{-N/2}
$$
for $t$ large enough. As $q<2$, this estimate and (\ref{lap2}) entail
that, for $t$ large enough,  
$$
\ell(t)\ge C\ t^{(4(2-q)-N q^2)/2q^2}
\mathop{\longrightarrow}_{t\to\infty} \infty\,, 
$$
since $q<4/(1+\sqrt{1+2N})$. Consequently, there exists $t_0$ large
enough such that $\ell(t_0)>K(q)$ and we may apply the first assertion
of Theorem \ref{th:05} to $t\longmapsto u(t_0+t)$ to complete the
proof. 
\qed

\bigskip

Under the assumptions of Theorem \ref{th:04}, part b) or Theorem \ref{th:05}, 
we may actually bound the $L^1$-norm of $u(t)$ from below and improve
significantly \cite[Proposition 2.1]{LS03}.   

\begin{prop} \label{pr:61}
Assume that $u_0$ satisfies \rf{u0-} and that 
$$
M_\infty := \lim_{t\to\infty} \Vert u(t)\Vert_\infty > 0\,.
$$
Then there is a constant $C=C(N,q,u_0)$ such that
\beqn
\label{volvic}
\Vert u(t)\Vert_1\ge C\ t^{N/q}\,,\quad t\ge 0\,.
\eeqn
\end{prop}

\noindent
{\sc Proof.} We fix $t>0$. For $k\ge 1$, let $x_k\in\RR^N$ be such
that $\Vert u(t)\Vert_\infty -1/k \le -u(x_k,t)$. For $R>0$, it
follows from \rf{grad} and the time monotonicity of $\Vert
u(t)\Vert_\infty$ that  
\bean
\Vert u(t)\Vert_1 & \ge & - \int_{\{ \vert x-x_k\vert\le R\}} u(x,t)\ dx \\
& \ge & \int_{\{ \vert x-x_k\vert\le R\}} \left( - u(x_k,t) - \vert
  x-x_k \vert\ \Vert\nabla u(t)\Vert_\infty \right)\ dx \\ 
& \ge & C\left( {R^N\over N}\ \left( \Vert u(t)\Vert_\infty -{1\over
      k} \right) - {C_1\ R^{N+1}\over N+1}\ \Vert
  u_0\Vert_\infty^{1/q}\ t^{-1/q} \right) \\ 
& \ge & C\ R^N\ \left( M_\infty - {1\over k} - C'\ R\ t^{-1/q}
\right)\,. 
\eean
Letting $k\to \infty$ and choosing $R= \left( M_\infty\ t^{1/q}
\right)/(2\ C')$ yields the claim \rf{volvic}. \qed 

%%%%%%%%%%%%%%%%%%%%%%%%%%%%%%%%%%%%%%%%%%%%%%%%%%%%%%%%%%%%%%%%%%%%%%%%%%
%%%%%
%%%%%%%%%%%%%%%%%%%%%%%%%%%%%%%%%%%%%%%%%%%%%%%%%%%%%%%%%%%%%%%%%%%%%%%%%%
%%%%%

\section{Proof of Theorem \ref{th:06} and Proposition \ref{pr:07}}
\setcounter{equation}{0}
\setcounter{remar}{0}

\noindent{\sc Proof of Theorem \ref{th:06}.}\\
\textbf{STEP 1.} 
Recall that, by \rf{u0-}, $u_0$ is a non-positive function.
We assume further that $u_0$ is compactly supported in a ball
$B(0,R_0)$ of $\RR^N$ for some $R_0>0$.  

For $\lambda\ge 1$, we introduce
$$
u_\lambda(x,t) := u(\lambda x, \lambda^q t)\,, \quad
(x,t)\in\RR^N\times (0,\infty)\,, 
$$
which solves
\beqn
\label{z5}
u_{\lambda,t} + \vert\nabla u_\lambda\vert^q = \lambda^{q-2}\ \Delta
u_\lambda \quad \mbox{ in }\quad \RR^N\times (0,\infty) 
\eeqn
with initial datum $u_\lambda(0)$. 

\begin{lem}\label{lez1}
There is a constant $C=C(N,q,\Vert u_0\Vert_\infty)$ such that, for
$t\ge 0$ and $\lambda\ge 1$,  
\beqn
\label{z6}
\Vert u_\lambda(t)\Vert_\infty + t^{1/q}\ \Vert \nabla
u_\lambda(t)\Vert_\infty + t\ \Vert u_{\lambda,t}(t)\Vert_\infty \le
C\,. 
\eeqn
\end{lem}

\noindent{\sc Proof.} It first follows from \rf{contrex} that 
$$
\Vert u_\lambda(t)\Vert_\infty = \Vert u(\lambda^q t)\Vert_\infty \le
\Vert u_0\Vert_\infty \,, 
$$
while Proposition \ref{prestgrad} yields 
$$
\Vert \nabla u_\lambda(t)\Vert_\infty = \lambda\ \Vert \nabla
u(\lambda^q t)\Vert_\infty \le C_1\ \Vert u_0\Vert_\infty^{1/q}\
t^{-1/q}\,. 
$$
We next infer from \cite[Theorem~5]{GGKpp} that 
$$
\Vert u_{\lambda,t}(t)\Vert_\infty = \lambda^q\ \Vert u_t(\lambda^q
t)\Vert_\infty \le \lambda^q\ C(N,q)\ \Vert u_0\Vert_\infty\ \left(
  \lambda^q t \right)^{-1} = C(N,q)\ \Vert u_0\Vert_\infty\ t^{-1}\,,  
$$
which completes the proof. \qed

\bigskip

Owing to Lemma \ref{lez1}, we may apply the Arzel\`a-Ascoli theorem
and deduce that there are a subsequence of $(u_\lambda)$ (not
relabeled) and a non-positive function $u_\infty\in
\mathcal{C}(\RR^N\times (0,\infty))$ such that  
\beqn
\label{z7}
u_\lambda\longrightarrow u_\infty \;\;\mbox{ in }\;\;
\mathcal{C}(B(0,R)\times (t_1,t_2)) 
\eeqn
for any $R>0$ and $0<t_1<t_2$. It also follows from (\ref{z7}) and
Lemma~\ref{lez1} that $u_\infty(t)\in \mathcal{BUC}(\RR^N)$ and
satisfies 
\beqn
\label{z8}
\Vert u_\infty(t)\Vert_\infty + t^{1/q}\ \Vert \nabla
u_\infty(t)\Vert_\infty + t\ \Vert u_{\infty,t}(t)\Vert_\infty \le C 
\eeqn
for each $t>0$. We next introduce the function $H_\lambda:
\RR\times\RR^N\times \mathcal{S}_N(\RR) \to \RR$  
defined by 
$$
H_\lambda(\xi_0,\xi,S) := \xi_0 + \vert\xi\vert^q - \lambda^{q-2}\
\mbox{tr}(S)\,, 
$$
where $\mathcal{S}_N(\RR)$ denotes the subset of symmetric matrices of
$\mathcal{M}_N(\RR)$ and $\mbox{tr}(S)$ denotes the trace of the
matrix $S$. On the one hand, we notice that \rf{z5} reads  
$$
H_\lambda(u_{\lambda,t}, \nabla u_\lambda, D^2 u_\lambda) = 0
\;\;\mbox{ in }\;\; \RR^N\times (0,\infty) 
$$
and that $H_\lambda$ is elliptic. On the other hand, $H_\lambda$
converges uniformly on every compact subset of $\RR\times\RR^N\times
\mathcal{S}_N(\RR)$ towards $H_\infty:\RR\times\RR^N\to \RR$ given by
$H_\infty(\xi_0,\xi) := \xi_0 + \vert\xi\vert^q$. Therefore, for every
$\tau>0$, $u_\infty(.+\tau)$ is the unique viscosity solution to
\rf{z4} with initial datum $u_\infty(\tau)$ ( see, e.g.,
\cite[Proposition~IV.1]{CL83} and \cite[Theorem~4.1]{CEL84}). In
addition, since $u_\infty(\tau)$ is bounded and Lipschitz continuous
by (\ref{z8}), we infer from \cite[Section~10.3, Theorem~3]{Ev98} that
$u_\infty(.+\tau)$ is given by the Hopf-Lax formula  
\beqn
\label{hlf}
u_\infty(x,t+\tau) = \inf_{y\in\RR^N} \left\{ u_\infty(y,\tau) +
  (q-1)\ q^{-q/(q-1)}\ \vert x-y\vert^{q/(q-1)}\ t^{-1/(q-1)} \right\} 
\eeqn 
for $(x,t)\in\RR^N\times [0,\infty)$.

It remains to identify the behavior of $u_\infty(t)$ as $t\to
0$. Consider first $x\in\RR^N$, $t\in (0,\infty)$ and $s\in (0,t)$. We
infer from \rf{lap2} and \rf{z5} that 
\bean
u_\lambda(x,t) & \le & u_\lambda(x,s) + \lambda^{q-2}\ \int_s^t \Delta
u_\lambda(x,\sigma)\ d\sigma \\ 
& \le & u_\lambda(x,s) + \lambda^{q-2}\ \int_s^t \lambda^2\  C\ \left(
  \lambda^q \sigma \right)^{-2/q}\ d\sigma \\ 
& \le & u_\lambda(x,s) - C\ \lambda^{q-2}\ \left( t^{(q-2)/q}
  -s^{(q-2)/q} \right)\,.
\eean
Since $q\in (1,2)$, we may pass to the limit as $\lambda\to\infty$ in
the previous inequality and use (\ref{z7}) to deduce that $t\mapsto
u_\infty(x,t)$ is non-increasing for every $x\in\RR^N$. Since
$u_\infty$ is bounded by (\ref{z8}), we may thus define $u_\infty(0)$
by  
\beqn
\label{zplus}
u_\infty(x,0) := \sup_{t>0}\, \{u_\infty(x,t)\} \in (-\infty,0]
\;\;\mbox{ for }\;\; x\in\RR^N\,. 
\eeqn
In particular, $u_\infty(0)$ is a lower semicontinuous function as the
supremum of continuous functions. 

More information on $u_\infty(0)$ are consequences of the next result. 

\begin{lem}\label{lez2}
For each $t>0$, there is $\varrho(t)>0$ such that $u_\infty(x,t)=0$ if
$\vert x\vert>\varrho(t)$ and 
\beqn
\label{z9}
\lim_{\lambda\to\infty} \Vert u_\lambda(t) - u_\infty(t)\Vert_\infty = 0\,.
\eeqn
Moreover, $\varrho(t)\to 0$ as $t\to 0$. 
\end{lem}

\medskip

Taking Lemma~\ref{lez2} for granted, we see that \rf{zplus} and
Lemma~\ref{lez2} imply that $u_\infty(x,0)=0$ for $x\ne 0$ since
$\varrho(t)\to 0$ as $t\to 0$. We set $\ell:=-u_\infty(0,0)$, so that 
$$
u_\infty(x,0) = - \ell \ \mathbf{1}_{\{0\}}(x)\,, \quad x\in\RR^N\,,
$$ and fix
$(x,t)\in\RR^N\times (0,\infty)$. We will now proceed along the lines
of \cite{St02} to show that $u_\infty(x,t) = Z_\ell(x,t)$ (recall that
$Z_\ell$ is defined in \rf{z3}). Introducing the notation $\mu:=
(q-1)\ q^{-q/(q-1)}\ t^{-1/(q-1)}$, it follows from (\ref{zplus}) and
Lemma~\ref{lez2} that, for $0<\sigma<\tau$ and $\vert
y\vert\le\varrho(\sigma)$,  
\bean
u_\infty(y,\sigma) + \mu\ \vert x-y\vert^{q/(q-1)} & \ge &
u_\infty(y,\tau) + \mu\ \vert x-y\vert^{q/(q-1)}\\ 
& \ge & u_\infty(0,\tau) + \mu\ \vert x\vert^{q/(q-1)} - \omega(\sigma)\,,
\eean
with
$$
\omega(\sigma) := \sup_{\vert y\vert\le\varrho(\sigma)} \left\vert
  u_\infty(y,\tau) - u_\infty(0,\tau) \right\vert + \mu\ \sup_{\vert
  y\vert\le\varrho(\sigma)} \left\vert \vert x-y\vert^{q/(q-1)} -
  \vert x\vert^{q/(q-1)} \right\vert\,, 
$$
while, for $0<\sigma<\tau$ and $\vert y\vert\ge\varrho(\sigma)$, 
$$
u_\infty(y,\sigma) + \mu\ \vert x-y\vert^{q/(q-1)} \ge 0\,.
$$
The previous bounds from below and (\ref{hlf}) entail that
$$
u_\infty(x,t+\sigma) \ge \min{\left\{ 0 , u_\infty(0,\tau) + \mu\
    \vert x\vert^{q/(q-1)} - \omega(\sigma) \right\}} 
$$
for $0<\sigma<\tau$. Since $\varrho(\sigma)\to 0$ as $\sigma\to 0$ and
$u_\infty\in\mathcal{C}(\RR^N\times (0,\infty))$, we may pass to the
limit as $\sigma\to 0$ in the above inequality and obtain 
$$
u_\infty(x,t) \ge \min{\left\{ 0 , u_\infty(0,\tau) + \mu\ \vert
    x\vert^{q/(q-1)} \right\}} 
$$
for $\tau>0$. Letting $\tau\to 0$ yields
$$
u_\infty(x,t) \ge \min{\left\{ 0 , -\ell + \mu\ \vert x\vert^{q/(q-1)}
  \right\}} = Z_\ell(x,t)\,. 
$$
On the other hand, (\ref{hlf}) and (\ref{zplus}) ensure that 
$$
u_\infty(x,t+\tau) \le \inf_{y\in\RR^N} \left\{ u_\infty(y,0) + \mu\
  \vert x-y\vert^{q/(q-1)} \right\} = Z_\ell(x,t)\,, 
$$
whence $u_\infty(x,t)\le Z_\ell(x,t)$ by the continuity of $u_\infty$
in $\RR^N\times (0,\infty)$. We have thus shown that
$u_\infty=Z_\ell$. In particular, $\Vert u_\infty(t)\Vert_\infty=\ell$
for $t\ge 0$. But \rf{z1} and \rf{z9} imply 
$$
\Vert u_\infty(t)\Vert_\infty = \lim_{\lambda\to\infty} \Vert
u_\lambda(t)\Vert_\infty = \lim_{\lambda\to\infty} \Vert u(\lambda^q
t)\Vert_\infty = M_\infty\,, 
$$
whence $\ell=M_\infty$ and $u_\infty=Z_{M_\infty}$. For $t>0$, the
sequence $(u_\lambda(t))$ has thus only one possible cluster point in
$L^\infty(\RR^N)$ as $\lambda\to\infty$, from which we conclude that
the whole family $(u_\lambda(t))$ converges to $Z_{M_\infty}(t)$ in
$L^\infty(\RR^N)$ as $\lambda\to\infty$. In particular, for $t=1$,  
$$
\lim_{\lambda\to\infty} \Vert u_\lambda(1) - Z_{M_\infty}(1)\Vert_\infty = 0\,.
$$
Setting $\lambda = t^{1/q}$ and using the self-similarity of
$Z_{M_\infty}$, we are finally led to \rf{z2}.  

\bigskip

\noindent\textbf{STEP 2.} We now consider an arbitrary function
$u_0\in\mathcal{C}_0(\RR^N)$ fulfilling \rf{u0-} and such that \rf{z1}
holds true. There is a sequence $(u_0^n)$ of non-positive functions in
$\mathcal{C}_c^\infty(\RR^N)$ such that 
$$
u_0^n \longrightarrow u_0 \;\;\mbox{ in }\;\; L^\infty(\RR^N)\,.
$$
For $n\ge 1$, we denote by $u^n$ the solution to (\ref{eq}) with
initial datum $u_0^n$ and put  
$$
M_\infty^n := \lim_{t\to\infty}\Vert u^n(t)\Vert_\infty\,.
$$
By \cite[Corollary~4.3]{GGK03}, we have
$$
\Vert u^n(t) - u(t) \Vert_\infty \le \Vert u_0^n - u_0\Vert_\infty
\;\;\mbox{ for }\;\; t\ge 0\,, 
$$
from which we readily deduce that 
$$
\left\vert M_\infty^n - M_\infty \right\vert \le \Vert u_0^n -
u_0\Vert_\infty\,. 
$$
Consequently, $M_\infty^n\longrightarrow M_\infty$ as $n\to \infty$
and (\ref{z1}) guarantees that $M_\infty^n>0$ for $n$ large
enough. The analysis performed in the previous step then implies that 
$$
\lim_{t\to\infty} \Vert u^n(t) - Z_{M_\infty^n}(t)\Vert_\infty=0
$$
for $n$ large enough. Therefore,
\bean
\Vert u(t) - Z_{M_\infty}(t)\Vert_\infty & \le & \Vert u(t) -
u^n(t)\Vert_\infty + \Vert u^n(t) - Z_{M_\infty^n}(t)\Vert_\infty \\ 
& + & \Vert Z_{M_\infty^n}(t) - Z_{M_\infty}(t)\Vert_\infty \\
& \le & \Vert u_0^n - u_0\Vert_\infty + \Vert u^n(t) -
Z_{M_\infty^n}(t)\Vert_\infty + \left\vert M_\infty^n - M_\infty
\right\vert \\ 
& \le & 2\ \Vert u_0^n - u_0\Vert_\infty + \Vert u^n(t) -
Z_{M_\infty^n}(t)\Vert_\infty\,, 
\eean
whence
$$
\limsup_{t\to\infty} \Vert u(t) - Z_{M_\infty}(t)\Vert_\infty \le 2\
\Vert u_0^n - u_0\Vert_\infty 
$$
for $n$ large enough. Letting $n\to\infty$ then completes the proof of
Theorem \ref{th:06}.  \hfill $\square$

\bigskip

\noindent{\sc Proof of Lemma~\ref{lez2}.} Let $h_0$ be a non-positive
function in $\mathcal{C}_c^\infty(\RR)$ such that $h_0(y) = -\Vert
u_0\Vert_\infty$ if $y\in (-R_0,R_0)$ (recall that $u_0$ is compactly
supported in $B(0,R_0)$). We denote by $h$ the solution to the
one-dimensional viscous Hamilton-Jacobi equation  
\bean
h_t - h_{yy} + \vert h_y\vert^q & = & 0 \;\;\mbox{ in }\;\; \RR\times
(0,\infty)\,,\\ 
h(0) & = & h_0 \;\;\mbox{ in }\;\; \RR\,.
\eean
For $i\in\{1,\ldots,N\}$ and $(x,t)\in\RR^N\times (0,\infty)$, we put
$h^i(x,t) := h(x_i,t)$ and notice that $h^i$ is the solution to
\rf{eq} with initial datum $h^i(0)\le u_0$. The comparison principle
then entails that  
\beqn
\label{z10}
h(x_i,t) = h^i(x,t) \le u(x,t) \le 0\,, \quad (x,t)\in \RR^N\times
(0,\infty)\,. 
\eeqn
We next introduce $w:=h_y$ and notice that $w$ is the solution to the
one-dimensional convection-diffusion equation  
\bear
\label{z11}
w_t - w_{yy} + \left( \vert w\vert^q \right)_y& = & 0 \;\;\mbox{ in
}\;\; \RR\times (0,\infty)\,,\\ 
\nonumber
w(0) & = & w_0 := h_{0,y} \;\;\mbox{ in }\;\; \RR\,.
\eear
The comparison principle then entails that 
\beqn
\label{z12}
b(y,t) \le w(y,t) \le a(y,t)\,, \quad (y,t)\in \RR\times (0,\infty)\,,
\eeqn
where $b\le 0$ and $a\ge 0$ denote the solutions to \rf{z11} with
initial data $b(0) = - w_0^-\le 0$ and $a(0)=w_0^+\ge 0$. Since
$w_0\in L^1(\RR)$, it follows from \cite{EVZ1} that 
\beqn
\label{z13}
\lim_{t\to\infty} \Vert b(t) - \Sigma_{-B}(t)\Vert_1 =
\lim_{t\to\infty} \Vert a(t) - \Sigma_{A}(t)\Vert_1 = 0\,,  
\eeqn
where $B:=\Vert b(0)\Vert_1$, $A:=\Vert a(0)\Vert_1$, and, for
$M\in\RR$, $\Sigma_M$ is the source solution to the one-dimensional
conservation law  
\bean
\Sigma_{M,t} + \left( \left\vert\Sigma_M\right\vert^q \right)_y & = &
0  \;\;\mbox{ in }\;\; \RR\times (0,\infty)\,,\\ 
\Sigma(0) & = & M\ \delta_0  \;\;\mbox{ in }\;\; \RR\,.
\eean
Here, $\delta_0$ denotes the Dirac mass in $\RR$ centered at
$y=0$. The source solution $\Sigma_M$ is actually given by  
$$
\Sigma_M(y,t) := y^{1/(q-1)}\  (qt)^{-1/(q-1)}\
\mathbf{1}_{[0,\xi_M(t)]}(y)\,, \quad \xi_M(t) := q\ \left(
  \frac{M}{q-1} \right)^{(q-1)/q}\ t^{1/q}\,,  
$$
if $M\ge 0$, and 
$$
\Sigma_M(y,t) := - \vert y\vert^{1/(q-1)}\  (qt)^{-1/(q-1)}\
\mathbf{1}_{[-\eta_M(t),0]}(y)\,, \quad \eta_M(t) := q\ \left(
  \frac{-M}{q-1} \right)^{(q-1)/q}\ t^{1/q}\,,  
$$
if $M\le 0$ (see, e.g., \cite{LP84}). In particular, $\Sigma_M$ satisfies
\beqn
\label{z14}
\lambda\ \Sigma_M(\lambda y,\lambda^q t) = \Sigma_M(y,t) \;\;\mbox{
  for }\;\; (\lambda,y,t)\in (0,\infty)\times\RR\times (0,\infty)\,. 
\eeqn

Now, let $t>0$ and set
$$
\varrho(t) := N^{1/2}\ \max{\left\{ \xi_A(t),\eta_{-B}(t)\right\}} \le
C\ t^{1/q}\,. 
$$
If $x\in\RR^N$ is such that $\vert x\vert>\varrho(t)$, there is
$i\in\{1,\ldots,N\}$ such that $\vert x_i\vert>\max{\left\{
    \xi_A(t),\eta_{-B}(t)\right\}}$, whence either $x_i>\xi_A(t)$ or
$x_i<-\eta_{-B}(t)$. In the latter case, we infer from \rf{z10},
\rf{z12} and \rf{z14} that 
\bean
0 \ge u_\lambda(x,t) & \ge & h(\lambda x_i,\lambda^q t) =
\int_{-\infty}^{\lambda x_i} w(y',\lambda^q t)\ dy' \\ 
& \ge & \lambda\ \int_{-\infty}^{x_i} b(\lambda y',\lambda^q t)\ dy'\\
& \ge & \lambda\ \int_{-\infty}^{x_i} \left( b(\lambda y',\lambda^q t)
  - \Sigma_{-B}(\lambda y',\lambda^q t) \right) \ dy'\\ 
& \ge & - \Vert (b-\Sigma_{-B})(\lambda^q t)\Vert_1\,.
\eean
Similarly, if $x_i>\xi_A(t)$, \rf{z10}, \rf{z12} and \rf{z14} yield
$$
0 \ge u_\lambda(x,t) \ge - \Vert (a-\Sigma_{A})(\lambda^q t)\Vert_1\,.
$$
Therefore, if $x\in\RR^N$ is such that $\vert x\vert>\varrho(t)$, then 
\beqn
\label{z15}
\vert u_\lambda(x,t) \vert \le \max{\left\{ \Vert
    (a-\Sigma_{A})(\lambda^q t)\Vert_1 , \Vert
    (b-\Sigma_{-B})(\lambda^q t)\Vert_1 \right\}}\,.  
\eeqn 
Passing to the limit as $\lambda\to\infty$ in \rf{z15} and using
    \rf{z7} and \rf{z13} provide the first assertion of Lemma
    \ref{lez2}. We next use once more \rf{z7} and \rf{z15} to conclude
    that \rf{z9} holds true. \hfill $\square$

\medskip

\noindent{\sc Proof of Proposition~\ref{pr:07}.} We keep the notations
of the proof of Theorem~\ref{th:06} and introduce
$$
U_\lambda(x,t) := u_{\lambda,x}(x,t) = \lambda\ u_x(\lambda x, \lambda^q
t)\,, \quad (x,t)\in\RR\times (0,\infty)\,.
$$
It follows from \rf{z5} and Lemma~\ref{lez1} that 
$$
U_{\lambda,t} + \left( \vert U_\lambda\vert^q \right)_x = \lambda^{q-2}\
U_{\lambda,xx}\,, \quad (x,t)\in\RR\times (0,\infty)\,,
$$
and
\beqn
\label{salsa}
\Vert U_\lambda(t)\Vert_1\le \Vert u_{0,x}\Vert_1 \;\;\mbox{ and }\;\;
t^{1/q}\ \Vert U_\lambda(t)\Vert_\infty \le C 
\eeqn 
for $t>0$. We recall that, by Theorem~\ref{th:06}, the family
$(u_\lambda)$ converges towards $Z_{M_\infty}$ in
$\mathcal{C}(\RR^N\times [t_1,t_2])$ for any $t_2>t_1>0$. Owing to
(\ref{salsa}), we readily conclude that $(U_\lambda)$
converges weakly-$\star$ towards $Z_{M_\infty,x}$ in
$L^\infty(\RR^N\times (t_1,t_2))$ for any $t_2>t_1>0$. We may then proceed
along the lines of \cite[Section~3]{EVZ1} to show that $(U_\lambda)$
converges towards $Z_{M_\infty,x}$ in $L^1(\RR)$ as
$\lambda\to\infty$. Expressing this convergence result in terms of
$U=u_x$ and using \rf{grad} yield Proposition~\ref{pr:07} by interpolation. 
\hfill $\square$

%%%%%%%%%%%%%%%%%%%%%%%%%%%%%%%%%%%%%%%%%%%%%%%%%%%%%%%%%%%%}
%%%%%%%%%%%%%%%%%%%%%%%%%%%%%%%%%%%%%%%%%%%%%%%%%%%%%%%%%%%%

\medskip

\noindent
{\bf Acknowledgements.}~ We thank Professor Herbert Koch for pointing
out to us Ref.~\cite{Ha93} and Professor Brian Gilding for useful comments on 
Proposition~\ref{prestlap}. The preparation of this paper was partially
supported by the KBN
grant 2 P03A 002 24, the POLONIUM project \'EGIDE--KBN No.~05643SE, and
the EU contract HYKE No.~HPRN-CT-2002-00282.

\end{document}